\documentclass[journal]{IEEEtran}

\usepackage{cite}
\ifCLASSINFOpdf
  \usepackage[pdftex]{graphicx}
  \usepackage{epstopdf}
  \usepackage{subfigure}
  \graphicspath{{./}{./figures/}}
\else
\fi

\usepackage{xcolor}

\usepackage{amsmath}
\usepackage{amssymb}

\def\R{\mathbb{R}}
\def\N{\mathbb{N}}

\renewcommand{\leq}{\leqslant}
\renewcommand{\geq}{\geqslant}

\newcommand*{\qed}{\hfill\ensuremath{\square}}%

\newtheorem{nnassumption}{\bf Assumption}

\newtheorem{nntheorem}{\bf Theorem}
\newenvironment{theorem}{\begin{nntheorem}\it}{\end{nntheorem}}
\newtheorem{nncorollary}{\bf Corollary}
\newenvironment{corollary}{\begin{nncorollary}\it}{\end{nncorollary}}
\newtheorem{nndefinition}{\bf Definition}

\newtheorem{nnproposition}{\bf Proposition}

\newtheorem{nnproblem}{\bf Problem}

\newtheorem{nnlemma}{\bf Lemma}
\newenvironment{lemma}{\begin{nnlemma}\it}{\end{nnlemma}}
\newtheorem{nnremark}{\bf Remark}
\newenvironment{remark}{\begin{nnremark} \rm }{\hfill \hspace*{1pt}\hfill $\circ$\end{nnremark}}

\allowdisplaybreaks

\begin{document}
%
\title{PI Regulation of a Reaction-Diffusion Equation with Delayed Boundary Control}
%
%
%

\author{Hugo~Lhachemi, Christophe~Prieur, Emmanuel~Tr{\'e}lat
\thanks{This publication has emanated from research supported in part by a research grant from Science Foundation Ireland (SFI) under grant number 16/RC/3872 and is co-funded under the European Regional Development Fund and by I-Form industry partners.
\newline\indent Hugo Lhachemi is with the School of Electrical and Electronic Engineering, University College Dublin, Dublin, Ireland (e-mail: hugo.lhachemi@ucd.ie).
\newline\indent Christophe Prieur is with Universit{\'e} Grenoble Alpes, CNRS, Grenoble-INP, GIPSA-lab, F-38000, Grenoble, France (e-mail: christophe.prieur@gipsa-lab.fr).
\newline\indent Emmanuel Tr{\'e}lat is with Sorbonne Universit\'e, CNRS, Universit\'e de Paris, Inria, Laboratoire Jacques-Louis Lions (LJLL), F-75005 Paris, France (email: emmanuel.trelat@sorbonne-universite.fr). 
}
}


%
%

\markboth{Manuscript}%
{Lhachemi \MakeLowercase{\textit{et al.}}: PI Regulation of a 1-D Linear Reaction-Diffusion Equation with Delayed Boundary Control}
%



\maketitle

\begin{abstract}
The general context of this work is the feedback control of an infinite-dimensional system so that the closed-loop system satisfies a fading-memory property and achieves the setpoint tracking of a given reference signal. More specifically, this paper is concerned with the Proportional Integral (PI) regulation control of the left Neumann trace of a one-dimensional reaction-diffusion equation with a delayed right Dirichlet boundary control. In this setting, the studied reaction-diffusion equation might be either open-loop stable or unstable. The proposed control strategy goes as follows. First, a finite-dimensional truncated model that captures the unstable dynamics of the original infinite-dimensional system is obtained via spectral decomposition. The truncated model is then augmented by an integral component on the tracking error of the left Neumann trace. After resorting to the Artstein transformation to handle the control input delay, the PI controller is designed by pole shifting. Stability of the resulting closed-loop infinite-dimensional system, consisting of the original reaction-diffusion equation with the PI controller, is then established thanks to an adequate Lyapunov function. In the case of a time-varying reference input and a time-varying distributed disturbance, our stability result takes the form of an exponential Input-to-State Stability (ISS) estimate with fading memory. Finally, another exponential ISS estimate with fading memory is established for the tracking performance of the reference signal by the system output. In particular, these results assess the setpoint regulation of the left Neumann trace in the presence of distributed perturbations that converge to a steady-state value and with a time-derivative that converges to zero. Numerical simulations are carried out to illustrate the efficiency of our control strategy.

\end{abstract}

\begin{IEEEkeywords}
1-D reaction-diffusion equation, PI regulation control, Neumann trace, Delay boundary control, Partial Differential Equations (PDEs).
\end{IEEEkeywords}

%
\IEEEpeerreviewmaketitle

\section{Introduction}\label{sec: Introduction}
%
%
%
%

\subsection{State of the art}

Motivated by the efficiency of Proportional-Integral (PI) controllers for the stabilization and regulation control of finite-dimensional systems, as well as its widespread adoption by  industry~\cite{aastrom1995pid,astrom2008feedback}, the opportunity of using PI controllers in the context of infinite-dimensional systems has attracted much attention in the recent years. One of the early attempts in this area was reported in~\cite{pohjolainen1982robust,pohjolainen1985robust}, then extended in~\cite{xu1995robust}, for bounded control operators. More recently, a number of works have been reported on the PI boundary control of linear hyperbolic systems~\cite{dos2008boundary,xu2014multivariable,bastin2015stability,lamare2015control}. The use of a PI boundary controller for 1-D nonlinear transport equation has been studied first in~\cite{trinh2017design} and then extended in~\cite{coron2019pi}. In particular, the former tackled the regulation problem for a constant reference input and in the presence of constant perturbations. The regulation of the downside angular velocity of a drilling string with a PI controller was reported in~\cite{terrand2018regulation}. The considered model consists of a wave equation coupled with ODEs in the presence of a constant disturbance. A related problem, the PI control of a drilling pipe under friction, was investigated in~\cite{barreau2019practical}. Recently, the opportunity to add an integral component to open-loop exponentially stable semigroups for the output tracking of a constant reference input and in the presence of a constant distributed perturbation was investigated in~\cite{terrand2018lyapunov,terrand2019adding} for unbounded control operators by using a Lyapunov functional design procedure. 

In this paper, we are concerned with the PI regulation control of the left Neumann trace of a one-dimensional reaction-diffusion equation with a delayed right Dirichlet boundary control. Specifically, we aim at achieving the setpoint reference tracking of a time-varying reference signal in spite of both the presence of an arbitrarily large constant input delay and a time-varying distributed disturbance. One of the early contributions regarding stabilization of PDEs with an arbitrarily large input delay deals with a reaction-diffusion equation~\cite{krstic2009control} where the controller was designed by resorting to the backstepping technique. A different approach, which is the one adopted in this paper, takes advantage of the following control design procedure initially reported in~\cite{russell1978controllability} and later used in~\cite{coron2004global,coron2006global,schmidt2006} to stabilize semilinear heat, wave or fluid equations via (undelayed) boundary feedback control: 1) design of the controller on a finite-dimensional model capturing the unstable modes of the original infinite-dimensional system; 2) use of an adequate Lyapunov function to assess that the designed control law stabilizes the whole infinite-dimensional system. The extension of this design procedure to the delay feedback control of a one-dimensional linear reaction-diffusion equation was reported in~\cite{prieur2018feedback}. The impact of the input-delay was handled in the control design by the synthesis of a predictor feedback via the classical Artstein transformation~\cite{artstein1982linear,richard2003time} (see also \cite{bresch2018new}). This control strategy was replicated in~\cite{guzman2019stabilization} for the feedback stabilization of a linear Kuramoto-Sivashinsky equation with delay boundary control. This idea was then generalized to the boundary feedback stabilization of a class of diagonal infinite-dimensional systems with delay boundary control for either a constant~\cite{lhachemi2018feedback,lhachemi2019control} or a time-varying~\cite{lhachemi2019lmi} input delay. 

\subsection{Investigated control problem}

Let $L > 0$, let $c \in L^{\infty}(0,L)$ and let $D > 0$ be arbitrary. We consider the one-dimensional reaction-diffusion equation over $(0,L)$ with delayed Dirichlet
boundary control:
\begin{subequations}
\begin{align}
& y_t = y_{xx} + c(x) y + d(t,x) , & (t,x) \in \mathbb{R}_+^* \times (0,L) \label{eq: investigated reaction-diffusion equation - beginning} \\
& y(t,0) = 0 , & t \geq 0 \\
& y(t,L) = u_D(t) \triangleq u(t-D) , & t \geq 0 \\
& y(0,x) = y_0(x) , & x \in (0,L) \label{eq: investigated reaction-diffusion equation - ending}
\end{align}
\end{subequations}
where $y(t,\cdot) \in L^2(0,L)$ is the state at time $t$, $u(t)\in\mathbb{R}$ is the control input, $D > 0$ is the (constant) control input delay, and $d(t,\cdot) \in L^2(0,L)$ is a time-varying distributed disturbance, continuously differentiable with respect to $t$.

In this paper, our objective is to achieve the PI regulation control of the left Neumann trace $y_x(t,0)$ to some prescribed reference signal, in the presence of the time-varying distributed disturbance $d$. More precisely, let $r : \mathbb{R}_+ \rightarrow \mathbb{R}$ be an arbitrary continuous function (reference signal). We aim at achieving the setpoint tracking of the time-varying reference signal $r(t)$ by the left Neumann trace $y_x(t,0)$. 

Note that an exponentially stabilizing controller for (\ref{eq: investigated reaction-diffusion equation - beginning}-\ref{eq: investigated reaction-diffusion equation - ending}) was designed in~\cite{prieur2018feedback} in the disturbance-free case ($d=0$) for a system trajectory evaluated in $H_0^1$-norm. The control strategy that we develop in the present paper elaborates on the one of \cite{prieur2018feedback}, adequately combined with a PI procedure. First, a finite-dimensional model capturing all unstable modes of the original infinite-dimensional system is obtained by an appropriate spectral decomposition. Following the standard PI approach, the tracking error on the left Neumann trace is then added as a new component to the resulting finite-dimensional system. Before synthetizing the PI controller, the control input delay is handled thanks to the Artstein transformation. A predictor feedback control, obtained by pole shifting, is then designed to exponentially stabilize the aforementioned truncated model. The core of the proof consists of establishing that this PI feedback controller exponentially stabilizes as well the complete infinite-dimensional system. This is done by an appropriate Lyapunov-based argument. 
The obtained results take the form of exponential Input-to-State Stability (ISS) estimates~\cite{sontag1989smooth} with fading memory of the reference input and the distributed perturbation. In the case where $r(t) \rightarrow r_e$, $d(t) \rightarrow d_e$ and $\dot{d}(t) \rightarrow 0$ when $t \rightarrow + \infty$, these estimates ensure the convergence of the state of the system, as well as the fulfillment of the desired setpoint regulation $y_x(t,0) \rightarrow r_e$. 

The paper is organized as follows. The proposed control strategy is introduced in Section~\ref{sec: control design strategy}. The study of the equilibrium points of the closed-loop system and the associated dynamics are presented in Section~\ref{sec: Equilibrium condition and related dynamics}. Then, the stability analysis of the closed-loop system is presented in Section~\ref{sec: stability analysis} while the assessment of the tracking performance is reported in Section~\ref{sec: setpoint reference tracking analysis}. The obtained results are illustrated by numerical simulations in Section~\ref{sec: numerical illustration}. Finally, concluding remarks are formulated in Section~\ref{sec: conclusion}.

\section{Control design strategy}\label{sec: control design strategy}

The sets of nonnegative integers, positive integers, real, nonnegative real, and positive real are denoted by $\mathbb{N}$, $\mathbb{N}^*$, $\mathbb{R}$, $\mathbb{R}_+$, and $\mathbb{R}_+^*$, respectively. All the finite-dimensional spaces $\mathbb{R}^p$ are endowed with the usual Euclidean inner product $\left<x,y\right> = x^\top y$ and the associated 2-norm $\Vert x \Vert = \sqrt{\left<x,x\right>} = \sqrt{x^\top x}$. For any matrix $M \in \mathbb{R}^{p \times q}$, $\Vert M \Vert$ stands for the induced norm of $M$ associated with the above 2-norms. For a given symmetric matrix $P \in \mathbb{R}^{p \times p}$, $\lambda_m(P)$ and $\lambda_M(P)$ denote its smallest and largest eigenvalues, respectively. In the sequel, the time derivative $\partial f/ \partial t$ is either denoted by $f_t$ or $\dot{f}$ while the spatial derivative $\partial f/ \partial x$ is either denoted by $f_x$ or $f'$.

\subsection{Augmented system for PI feedback control}
The control design objective is: 1) to stabilize the reaction-diffusion system (\ref{eq: investigated reaction-diffusion equation - beginning}-\ref{eq: investigated reaction-diffusion equation - ending}); 2) to ensure the setpoint tracking of the reference signal $r(t)$ by the left Neumann trace $y_x(t,0)$. We address this problem by designing a PI controller. Following the general PI scheme, we augment the system by introducing a new state $z(t)\in\R$ taking the form of the integral of the tracking error $y_x(t,0) - r(t)$ (as for finite-dimensional systems, the objective of this integral component is to ensure the setpoint tracking of the reference signal in the presence of the distributed disturbance $d$):
\begin{subequations}
\begin{align}
& y_t = y_{xx} + c(x) y + d(t,x) , & (t,x) \in \mathbb{R}_+^* \times (0,L) \label{eq: augmented reaction-diffusion equation - beginning} \\
& \dot{z}(t) = y_x(t,0) - r(t), & t \geq 0 \\
& y(t,0) = 0 , & t \geq 0 \\
& y(t,L) = u_D(t) \triangleq u(t-D) , & t \geq 0 \\
& y(0,x) = y_0(x) , & x \in (0,L) \\
& z(0) = z_0   \label{eq: augmented reaction-diffusion equation - ending}
\end{align}
\end{subequations}
where $z_0 \in \mathbb{R}$ stands for the initial condition of the integral component. As we are only concerned in prescribing the future of the system, we assume that the system is uncontrolled for $t < 0$, i.e., $u(t) = 0$ for $t < 0$. Consequently, due to the input delay $D > 0$, the system is in open loop over the time range $[0,D)$ as the impact of the control strategy actually applies in the boundary condition only for $t \geq D$.

\subsection{Modal decomposition}
It is convenient to rewrite (\ref{eq: augmented reaction-diffusion equation - beginning}-\ref{eq: augmented reaction-diffusion equation - ending}) as an equivalent homogeneous Dirichlet problem. Specifically, assuming\footnote{This property will be ensured by the construction carried out in the sequel.} that $u$ is continuously differentiable and setting $w(t,x) = y(t,x) - \frac{x}{L} u_D(t)$, we have
\begin{subequations}
\begin{align}
& w_t = w_{xx} + c(x) w + \frac{x}{L} c(x) u_D - \frac{x}{L} \dot{u}_D(t) + d(t,x) \label{eq: homogeneous reaction-diffusion equation - PDE} \\
& \dot{z}(t) = w_x(t,0) + \frac{1}{L} u_D(t) - r(t) \label{eq: homogeneous reaction-diffusion equation - integral component} \\
& w(t,0) = w(t,L) = 0 \label{eq: homogeneous reaction-diffusion equation - BC} \\
& w(0,x) = y_0(x) - \frac{x}{L} u_D(0) \label{eq: homogeneous reaction-diffusion equation - IC PDE} \\
& z(0) = z_0 \label{eq: homogeneous reaction-diffusion equation - IC integral component}
\end{align}
\end{subequations}
for $t > 0$ and $x \in (0,1)$. 
We consider the real state-space $L^2(0,1)$ endowed with its usual inner product $\left\langle f , g \right\rangle = \int_0^L f(x) g(x) \,\mathrm{d}x$. Introducing the operator $\mathcal{A} = \partial_{xx} + c\, \mathrm{id} : D(\mathcal{A}) \subset L^2(0,L) \rightarrow L^2(0,L)$ defined on the domain $D(\mathcal{A}) = H^2(0,L) \cap H_0^1(0,L)$, (\ref{eq: homogeneous reaction-diffusion equation - PDE}-\ref{eq: homogeneous reaction-diffusion equation - BC}) can be rewritten as
\begin{subequations}
\begin{align}
& w_t(t,\cdot) = \mathcal{A} w(t,\cdot)+a(\cdot)u_D(t)+b(\cdot)\dot{u}_D(t)+d(t,\cdot) \label{eq: abstract form - 1} \\
& \dot z(t) = w_x(t,0)+\frac{1}{L}u_D(t)-r(t) \label{eq: abstract form - 2}
\end{align}
\end{subequations}
with $a(x)=\frac{x}{L}c(x)$ and $b(x)=-\frac{x}{L}$ for every $x\in(0,L)$, with initial conditions (\ref{eq: homogeneous reaction-diffusion equation - IC PDE}-\ref{eq: homogeneous reaction-diffusion equation - IC integral component}). Since $\mathcal{A}$ is self-adjoint and of compact resolvent, we consider a Hilbert basis $(e_j)_{j\geq 1}$ of $L^2(0,L)$ consisting of eigenfunctions of $\mathcal{A}$ associated with the sequence of real eigenvalues
\begin{equation*}
-\infty<\cdots<\lambda_j<\cdots<\lambda_1\quad \textrm{with}\quad \lambda_j\underset{j\rightarrow +\infty}{\longrightarrow}-\infty .
\end{equation*}
Note that $e_j(\cdot)\in H^1_0(0,L)\cap C^2([0,L])$ for every $j\geq 1$ and
\begin{equation}\label{eq: symptotic behaviors}
e_j'(0) \sim \sqrt{\frac{2}{L}} \sqrt{\vert \lambda_j \vert} ,
\quad \lambda_j \sim -\frac{\pi^2 j^2}{L^2} ,
\end{equation}
when $j \rightarrow + \infty$.
The solution $w(t,\cdot)\in H^2(0,L)\cap H^1_0(0,L)$ of (\ref{eq: abstract form - 1}) can be expanded as a series in the eigenfunctions $e_j(\cdot)$, convergent in $H_0^1(0,L)$,
\begin{equation}\label{eq: series expansion solutions}
w(t,\cdot)=\sum_{j=1}^{+\infty}w_j(t)e_j(\cdot) .
\end{equation}
Therefore (\ref{eq: abstract form - 1}-\ref{eq: abstract form - 2}) is equivalent to the infinite-dimensional control system:
\begin{subequations}
\begin{align}
\dot{w}_j(t) & = \lambda_j w_j(t) + a_j u_D(t) + b_j \dot{u}_D(t) + d_j(t) \label{eq: prel 1 spectral decomposition - 1} \\
\dot{z}(t) & = \sum_{j \geq 1} w_j(t) e_j'(0) + \frac{1}{L} u_D(t) - r(t) \label{eq: prel 1 spectral decomposition - 2}
\end{align}
\end{subequations}
for $j\in\N^*$, with
\begin{align*}
w_j(t) & = \left\langle w(t,\cdot) , e_j \right\rangle = \int_0^L w(t,x) e_j(x) \,\mathrm{d}x , \\
a_j & = \left\langle a , e_j \right\rangle = \frac{1}{L}\int_0^L x c(x) e_j(x) \,\mathrm{d}x , \\
b_j & = \langle b , e_j \rangle = -\frac{1}{L} \int_0^L x e_j(x) \,\mathrm{d}x , \\
d_j(t) & = \langle d(t,\cdot) , e_j \rangle = \int_0^L d(t,x) e_j(x) \,\mathrm{d}x .
\end{align*}
Introducing the auxiliary control input $v \triangleq \dot{u}$, and denoting $v_D(t) \triangleq v(t-D)$, (\ref{eq: prel 1 spectral decomposition - 1}-\ref{eq: prel 1 spectral decomposition - 2}) can be rewritten as
\begin{subequations}
\begin{align}
\dot{u}_D(t) & = v_D(t) \label{eq: prel 2 spectral decomposition - 1} \\
\dot{w}_j(t) & = \lambda_j w_j(t) + a_j u_D(t) + b_j v_D(t) + d_j(t) \label{eq: prel 2 spectral decomposition - 2} \\
\dot{z}(t) & = \sum_{j \geq 1} w_j(t) e_j'(0) + \frac{1}{L} u_D(t) - r(t) \label{eq: prel 2 spectral decomposition - 3}
\end{align}
\end{subequations}
for $j \in\N^*$. Since $u(t) = 0$ for $t < 0$, (\ref{eq: prel 2 spectral decomposition - 1}) imposes that the auxiliary control input is such that $v(t) = 0$ for $t < 0$, and that the corresponding initial condition satisfies $u_D(0)=u(-D)=0$. In the sequel, we design the control law $v$ in order to stabilize (\ref{eq: prel 2 spectral decomposition - 1}-\ref{eq: prel 2 spectral decomposition - 3}). In this context, the actual control input $u$ associated with the original system (\ref{eq: augmented reaction-diffusion equation - beginning}-\ref{eq: augmented reaction-diffusion equation - ending}) is $u(t) = \int_0^t v(\tau) \,\mathrm{d}\tau$ for every $t \geq 0$.

\subsection{Finite-dimensional truncated model}
In what follows, we fix the integer $n \in\mathbb{N}$ such that $\lambda_{n+1} < 0 \leq \lambda_n$. In particular, we have $\lambda_{j} \geq 0$ when $1 \leq j \leq n$ and $\lambda_{j} \leq \lambda_{n+1} < 0$ when $j \geq n+1$. 

\begin{remark}
In the case of an open-loop stable reaction-diffusion equation, we have $n=0$. In this particular case, as discussed in the sequel, the objective of the control design is to ensure the output regulation while preserving the stability of the closed-loop system.
\end{remark}

Let us first show how to obtain a finite-dimensional truncated model capturing the $n$ first modes of the reaction diffusion-equation. We follow \cite{prieur2018feedback}. Setting
\begin{equation*}
X_1(t) = \begin{pmatrix} u_D(t) \\ w_1(t) \\ \vdots \\ w_n(t)
\end{pmatrix} ,
\quad 
A_1 = \begin{pmatrix}
0      &         0 & \cdots &         0 \\
a_1    & \lambda_1 & \cdots &         0 \\
\vdots &    \vdots & \ddots &    \vdots \\
a_n    &         0 & \cdots & \lambda_n 
\end{pmatrix} ,
\end{equation*}
\begin{align*}
B_1 & = \begin{pmatrix} 1 & b_1 & \ldots & b_n \end{pmatrix}^\top ,
\\
D_1(t) & = \begin{pmatrix} 0 & d_1(t) & \ldots & d_n(t) \end{pmatrix}^\top , 
\end{align*}
with $X_1(t) \in \mathbb{R}^{n+1}$, $A_1 \in \mathbb{R}^{(n+1) \times (n+1)}$, $B_1 \in \mathbb{R}^{n+1}$, $D_1(t) \in \mathbb{R}^{n+1}$, (\ref{eq: prel 2 spectral decomposition - 1}) and the $n$ first equations of (\ref{eq: prel 2 spectral decomposition - 2}) yield
\begin{equation}\label{eq: ODE prel truncated model}
\dot{X}_1(t) = A_1 X_1(t) + B_1 v_D(t) + D_1(t) .
\end{equation}
We could now augment the state-vector $X_1$ to include the integral component $z$ in the control design. However, the time derivative of $z$, given by (\ref{eq: prel 2 spectral decomposition - 3}), involves all coefficients $w_j(t)$, $j \geq 1$. Thus, the direct augmentation of the state vector $X_1$ with the integral component $z$ does not allow the derivation of an ODE involving only the $n$ first modes of the reaction-diffusion equation. To overcome this issue, we set
\begin{equation}\label{series10}
\zeta(t) \triangleq z(t) - \sum_{j \geq n+1} \frac{e_j'(0)}{\lambda_j} w_j(t) .
\end{equation}
Noting that $\left\vert \frac{e_j'(0)}{\lambda_j} \right\vert^2 \sim \frac{2L}{\pi^2 j^2}$ when $j \rightarrow + \infty$ and thus that $(e'_j(0)/\lambda_j)_j$ and $(w_j(t))_j$ are square summable sequences, using the Cauchy-Schwarz inequality, we see that the series \eqref{series10} is convergent and that
\begin{align*}
\dot{\zeta}(t) 
& = \dot{z}(t) - \sum_{j \geq n+1} \frac{e_j'(0)}{\lambda_j} \dot{w}_j(t) \\
& = \alpha u_D(t) + \beta v_D(t) - \gamma(t) +\sum_{j=1}^n w_j(t)e_j'(0) , 
\end{align*}
where we have used (\ref{eq: prel 2 spectral decomposition - 2}-\ref{eq: prel 2 spectral decomposition - 3}), with
\begin{subequations}
\begin{equation}\label{eq: def parameters dot_zeta - 1}
\alpha = \frac{1}{L} - \sum_{j \geq n+1} \frac{e_j'(0)}{\lambda_j} a_j, \quad
\beta = - \sum_{j \geq n+1} \frac{e_j'(0)}{\lambda_j} b_j, \quad
\end{equation}
\begin{equation}\label{eq: def parameters dot_zeta - 2}
\gamma(t) = r(t) + \sum_{j \geq n+1} \frac{e_j'(0)}{\lambda_j} d_j(t) .
\end{equation}
\end{subequations}
The convergence of the above series follow again by the Cauchy-Schwarz inequality. Then we have
\begin{equation}\label{eq: zeta prel truncated model}
\dot{\zeta}(t) = L_1 X_1(t) + \beta v_D(t) - \gamma(t) 
\end{equation}
with $L_1 = \begin{pmatrix} \alpha & e'_1(0) & \ldots & e'_{n}(0) \end{pmatrix} \in \mathbb{R}^{1 \times (n+1)}$.
Now, defining the augmented state-vector $X(t) = \begin{bmatrix} X_1(t)^\top & \zeta(t) \end{bmatrix}^\top \in \mathbb{R}^{n+2}$, the exogenous input $\Gamma(t) = \begin{bmatrix} D_1(t)^\top & -\gamma(t) \end{bmatrix}^\top \in \mathbb{R}^{n+2} $ and the matrices
\begin{equation}\label{eq: ODE truncated model - matrices A and B}
A = 
\begin{pmatrix}
A_1 & 0 \\
L_1 & 0
\end{pmatrix}
\in \mathbb{R}^{(n+2) \times (n+2)}
, \quad
B = 
\begin{pmatrix}
B_1\\ \beta
\end{pmatrix}
\in \mathbb{R}^{n+2} ,
\end{equation} 
we obtain from (\ref{eq: ODE prel truncated model}) and (\ref{eq: zeta prel truncated model}) the control system
\begin{equation}\label{eq: ODE truncated model}
\dot{X}(t) = A X(t) + B v_D(t) + \Gamma(t) 
\end{equation}
which is the finite-dimensional truncated model capturing the unstable part of the infinite-dimensional augmented with an integral component for generating the actual control input and an integral component for setpoint reference tracking. In particular, the system (\ref{eq: ODE truncated model}) only involves the $n$ first modes of the reaction-diffusion equation.

\begin{remark}
The above developments allow the particular case $n = 0$, which corresponds to the configuration where (\ref{eq: investigated reaction-diffusion equation - beginning}-\ref{eq: investigated reaction-diffusion equation - ending}) is open-loop stable. In this configuration, the vectors and matrices of the truncated model (\ref{eq: ODE truncated model}) reduce to $X(t) = \begin{bmatrix} u_D(t) & \zeta(t) \end{bmatrix}^\top \in \mathbb{R}^2$, $\Gamma(t) = \begin{bmatrix} 1 & - \gamma(t) \end{bmatrix}^\top \in \mathbb{R}^2$,
\begin{equation*}
A
= \begin{pmatrix}
0 & 0 \\ \alpha & 0
\end{pmatrix} ,
\quad
B
= \begin{pmatrix}
1 \\ \beta
\end{pmatrix} .
\end{equation*} 
In this setting, the control objective consists of ensuring the setpoint tracking of the system output $y_x(t,0)$ while preserving the stability of the closed-loop system.
\end{remark}

Putting together the finite-dimensional truncated model (\ref{eq: ODE truncated model}) along with (\ref{eq: prel 2 spectral decomposition - 2}) for $j \geq n+1$ which correspond to the modes of the original infinite-dimensional system neglected by the truncated model, we get the final representation used for both control design and stability analyses:
\begin{subequations}
\begin{align}
\dot{X}(t) & = A X(t) + B v_D(t) + \Gamma(t) \label{eq: final model - 1} \\
\dot{w}_j(t) & = \lambda_j w_j(t) + a_j u_D(t) + b_j v_D(t) + d_j(t) \label{eq: final model - 2} 
\end{align}
\end{subequations}
with $j\geq n+1$.

\subsection{Controllability of the finite-dimensional truncated model}
As mentioned in the introduction, the control design strategy relies now on the two following steps. First, we want to design a controller for the finite-dimensional system (\ref{eq: ODE truncated model}). Second, we aim at assessing that the obtained PI controller successfully stabilizes the original infinite-dimensional system (\ref{eq: augmented reaction-diffusion equation - beginning}-\ref{eq: augmented reaction-diffusion equation - ending}) and provides the desired setpoint reference tracking. In order to fulfill the first objective, we first establish the controllability property for the pair $(A,B)$.

\begin{lemma}
The pair $(A,B)$ satisfies the Kalman condition.
\end{lemma}

\textbf{Proof:}
Considering the structures of $A$ and $B$ defined by (\ref{eq: ODE truncated model - matrices A and B}), we apply Lemma~\ref{lem: technical lemma} reported in Appendix. More specifically, from the implication $(i) \Rightarrow (ii)$, we need to check that the pair $(A_1,B_1)$ satisfies the Kalman condition and the square matrix $\begin{pmatrix} A_1 & B_1 \\ L_1 & \beta \end{pmatrix}$ is invertible. The first condition is indeed true as straightforward computations show that $\mathrm{det}\left( B_1,A_1B_1,\ldots,A_1^nB_1\right)=\prod_{j=1}^n (a_j+\lambda_jb_j) \mathrm{VdM}(\lambda_1,\ldots,\lambda_n) \neq 0$,
where $\mathrm{VdM}$ is a Vandermonde determinant, because all eigenvalues are distinct and, using $\mathcal{A} e_j = \lambda_j e_j$ and an integration by parts, $a_j + \lambda_j b_j = -e_j'(L) \neq 0$ by Cauchy uniqueness (see also \cite{prieur2018feedback}). Thus, we focus on the invertibility condition:
\begin{align*}
\mathrm{det}\begin{pmatrix}
A_1 & B_1\\
L_1 & \beta
\end{pmatrix}
& =
\mathrm{det}
\begin{pmatrix}
0         &       0         & \cdots &    0  & 1         \\
a_1 & \lambda_1 & \cdots &    0    & b_1       \\
\vdots    &  \vdots         & \ddots &   \vdots   & \vdots    \\
a_n &  0              & \cdots & \lambda_n & b_n \\
\alpha & e_1'(0) & \cdots & e_n'(0) & \beta
\end{pmatrix} \\
& = (-1)^{n+1} \mathrm{det}
\begin{pmatrix}
a_1 & \lambda_1 & \cdots &    0           \\
\vdots    &  \vdots         & \ddots &   \vdots      \\
a_n &  0              & \cdots & \lambda_n  \\
\alpha & e_1'(0) & \cdots & e_n'(0) 
\end{pmatrix} .
\end{align*}
We now consider two distinct cases depending on whether $\lambda = 0$ is an eigenvalue of $\mathcal{A}$ or not.

Let us first consider the case where $\lambda = 0$ is not an eigenvalue of $\mathcal{A}$. In particular, $\lambda_1,\ldots,\lambda_n$ are all non zero and thus row operations applied to the last row yield:
\begin{align*}
& \mathrm{det}\begin{pmatrix}
A_1 & B_1\\
L_1 & \beta
\end{pmatrix} \\[-6mm]
& \qquad =
(-1)^{n+1}
\mathrm{det}
\begin{pmatrix}
a_1 & \lambda_1 & \cdots &    0           \\
\vdots    &  \vdots         & \ddots &   \vdots      \\
a_n &  0              & \cdots & \lambda_n  \\
\alpha - \sum\limits_{i=1}^n a_i \frac{e_i'(0)}{\lambda_i} & 0 & \cdots & 0 
\end{pmatrix} \\
& \qquad = - \left( \alpha - \sum\limits_{i=1}^n a_i \frac{e_i'(0)}{\lambda_i} \right) \prod\limits_{j=1}^n \lambda_j .
\end{align*}
Consequently, based on the definition of the constant $\alpha$ given by (\ref{eq: def parameters dot_zeta - 1}), the above determinant is not zero if and only if
\begin{equation}\label{eq: commandability - first case - condition}
\sum\limits_{j \geq 1} a_i \frac{e_i'(0)}{\lambda_i} \neq \frac{1}{L} .
\end{equation}
We note that this condition is independent of the number $n$ of modes of the infinite-dimensional system captured by the truncated model and we show in the sequel that (\ref{eq: commandability - first case - condition}) always holds true. To do so, let $y_e$ be the stationary solution of (\ref{eq: investigated reaction-diffusion equation - beginning}-\ref{eq: investigated reaction-diffusion equation - ending}) associated with the constant boundary input $u_e = 1$ and zero distributed disturbance, i.e., $(y_e)_{xx} + c y_e =0$ with $y_e(0) = 0$ and $y_e(L)=1$. Such a function $y_e$ indeed exists and can be obtained as follows. By assumption, $\lambda = 0$ is not an eigenvalue of $\mathcal{A}$. Thus the solution $y_0$ of $(y_0)_{xx} + c y_0 =0$ with $y_0(0) = 0$ and $y'_0(0) = 1$ satisfies $y_0(L) \neq 0$. Hence, one can obtain the claimed function by defining $y_e(x) = y_0(x) / y_0(L)$. Now, $w_e(x) \triangleq y_e(x) - \frac{x}{L}$ is a stationary solution of (\ref{eq: homogeneous reaction-diffusion equation - PDE}) and (\ref{eq: homogeneous reaction-diffusion equation - BC}-\ref{eq: homogeneous reaction-diffusion equation - IC PDE}) in the sense that $(w_e)_{xx} + c w_e + \frac{x}{L} c = 0$ with $w_e(0) = w_e(L) = 0$. From (\ref{eq: prel 1 spectral decomposition - 1}), $\lambda_j w_{e,j} + a_j = 0$ and thus $w_{e,j} = - \frac{a_j}{\lambda_j}$. We deduce that
\begin{equation*}
(w_{e})_x(0) = \sum\limits_{j \geq 1} w_{e,j} e'_j(0) = - \sum\limits_{j \geq 1} \frac{a_j}{\lambda_j} e'_j(0) .
\end{equation*}
Hence (\ref{eq: commandability - first case - condition}) holds if and only if $(w_{e})_x(0) \neq - \frac{1}{L}$, which is equivalent to $(y_{e})_x(0) \neq 0$. By Cauchy uniqueness, the condition $(y_{e})_x(0) = 0$, along with  $(y_e)_{xx} + c y_e =0$ and $y_e(0) = 0$, implies that $y_e = 0$, which contradicts $y_e(L)=1$. Thus (\ref{eq: commandability - first case - condition}) holds and the system is controllable.

Let us now consider the second case, i.e., $\lambda = 0$ is an eigenvalue of $\mathcal{A}$. Based on the definition of the integer $n$, we have $n \geq 1$ and $\lambda_n = 0$ while $\lambda_k > 0$ for all $1 \leq k \leq n-1$. Expanding the determinant, first, along the $(n+1)$-th column, and then, along the $n$-th row, we obtain
\begin{equation*}
\mathrm{det}\begin{pmatrix}
A_1 & B_1\\
L_1 & \beta
\end{pmatrix}
= \,a_n\, e_n'(0)\, \prod_{i=1}^{n-1} \lambda_i .
\end{equation*}
By Cauchy uniqueness, we have $e_n'(0)\neq 0$ (otherwise $e_n$ would be solution of a second-order ODE with the boundary conditions $e_n(0)=e'_n(0)=0$, yielding the contradiction $e_n = 0$). Thus, the above determinant is nonzero if and only if $a_n \neq 0$. We proceed by contradiction. Using $e''_n + c e_n = \mathcal{A} e_n = 0$ and $a_n = \frac{1}{L} \int_0^L x c(x) e_n(x) \,\mathrm{d}x = 0$, we obtain by integration by parts: $0 = - \int_0^L x e''_n(x) \,\mathrm{d}x = - \left[ x e'_n(x) \right]_{x=0}^{x=L} + \int_{0}^{L} e'_n(x) \,\mathrm{d}x = - L e'_n(L)$, whence $e'_n(L) = 0$. This result, along with $e''_n + c e_n = 0$ and $e_n(L) = 0$, yields by Cauchy uniqueness the contradiction $e_n = 0$. Thus, $a_n \neq 0$ and the system is controllable.
\qed

\subsection{Control design strategy}
Using the controllability property of the pair $(A,B)$, we propose to resort to the classical predictor feedback to stabilize the finite-dimensional truncated model (\ref{eq: final model - 1}). Specifically, introducing the Artstein transformation
\begin{equation}\label{eq: Artstein transformation}
Z(t) = X(t) + \int_{t-D}^{t} e^{A (t-D-\tau)} B v(\tau) \,\mathrm{d}\tau
\end{equation}
(see \cite{artstein1982linear}), 
straightforward computations show that
\begin{equation*}
\dot{Z}(t) = A Z(t) + e^{-D A} B v(t) + \Gamma(t) .
\end{equation*} 
Since $(A,B)$ satisfies the Kalman condition, the pair $(A,e^{-D A}B)$ also satisfies the Kalman condition and we infer the existence of a feedback gain $K \in \mathbb{R}^{1 \times (n+2)}$ such that $A_K \triangleq A + e^{-D A} B K$ is Hurwitz. We choose the control law
\begin{equation}\label{eq: control input}
v(t) = \chi_{[0,+\infty)}(t) K Z(t) 
\end{equation}
where $\chi_{[0,+\infty)}$ denotes the characteristic function of the interval $[0,+\infty)$, which is used to capture the fact that we are only concerned by imposing a non zero control input for $t > 0$. Then  we obtain the stable closed-loop dynamics
\begin{equation*}
\dot{Z}(t) = A_K Z(t) + \Gamma(t) .
\end{equation*} 

\begin{remark}
The first component of $Z(t)$ is $u(t)$. Indeed, denoting by $E_1 = \begin{bmatrix} 1 & 0 & \ldots & 0 \end{bmatrix} \in \mathbb{R}^{1 \times (n+2)}$, we have
\begin{align*}
E_1 Z(t)
& = E_1 X(t)
+ \int_{t-D}^t E_1 e^{(t-s-D)A} B v(s) \,\mathrm{d}s \\
& = u_D(t) + \int_{t-D}^t v(s) \,\mathrm{d}s 
 = u(t) 
\end{align*}
where we have used that the first row of $A$ is null, that $v(t) = \dot{u}(t)$ for $t \geq 0$ and that $u(t) = 0$ when $t \leq 0$.
\end{remark}

\begin{remark}
Putting together (\ref{eq: Artstein transformation}-\ref{eq: control input}) and using the fact that $v(t) = 0$ for $t \leq 0$, we obtain that the control input $v$ is solution of the fixed point implicit equation
\begin{equation*}
v(t) = \chi_{[0,+\infty)}(t) K X(t) + K \int_{\max(t-D,0)}^{t}\!\!\!\! e^{A (t-D-\tau)} B v(\tau) \,\mathrm{d}\tau .
\end{equation*}
Existence and uniqueness of the solution of the above equation as well as regularity properties and  inversion of the Artstein transformation are reported in~\cite{bresch2018new}.
\end{remark}

The main objective is now to establish that the feedback control (\ref{eq: control input}) stabilizes as well the original infinite-dimensional system (or, in the case $n = 0$, preserves the stability property of the system) while providing a setpoint tracking of the time-varying reference signal $r(t)$ by the left Neumann trace $y_x(t,0)$. The former is studied in Section~\ref{sec: stability analysis} while the latter is investigated in Section~\ref{sec: setpoint reference tracking analysis}.

\section{Equilibrium condition and related dynamics}\label{sec: Equilibrium condition and related dynamics}

In the sequel, $r_e \in\R$ and $d_e \in L^2(0,L)$ stand for ``nominal values'' of the time-varying reference signals $r(t)$ and the distributed disturbance $d(t)$, respectively. Even if $r_e$ and $d_e$ can be selected arbitrarily, the two following (distinct) cases will be of particular interest in the sequel:
\begin{itemize}
\item $\vert r(t) - r_e \vert \leq \delta_r$ and $\Vert d(t) - d_e \Vert \leq \delta_d$ for some  $\delta_r , \delta_d > 0$;
\item $r(t) \rightarrow r_e$ and $d(t) \rightarrow d_e$ when $t \rightarrow + \infty$.
\end{itemize}

\subsection{Characterization of equilibrium for the closed-loop system}
Setting $d_{e,j} = \left\langle d_e , e_j \right\rangle = \int_0^L d_e(x) e_j(x) \,\mathrm{d}x$ for $j \geq 1$, 
$\Delta r = r - r_e$, $\Delta d = d - d_e$, $\Delta d_j = d_j - d_{e,j}$, 
{\small
\begin{equation*}
\Gamma_e 
= \begin{pmatrix} 0 \\ d_{e,1} \\ \vdots \\ d_{e,n} \\ - r_e - \sum\limits_{j \geq n+1} \frac{e_j'(0)}{\lambda_j} d_{e,j} \end{pmatrix} ,
\,
\Delta\Gamma 
= \begin{pmatrix} 0 \\ \Delta d_1 \\ \vdots \\ \Delta d_n \\ - \Delta r - \sum\limits_{j \geq n+1} \frac{e_j'(0)}{\lambda_j} \Delta d_j \end{pmatrix}
\end{equation*}
}
we obtain from (\ref{eq: final model - 1}-\ref{eq: final model - 2}) and (\ref{eq: control input})
\begin{align*}
\dot{Z}(t) & = A_K Z(t) + \Gamma_e + \Delta\Gamma(t) \\
\dot{w}_j(t) & = \lambda_j w_j(t) + a_j u_D(t) + b_j v_D(t) + d_{e,j} + \Delta d_j(t)
\end{align*}
for $j \geq n+1$. We now characterize the equilibrium condition of the above closed-loop system associated with the constant reference input $r(t) = r_e \in \mathbb{R}$ and the constant distributed disturbance $d(t) = d_e \in L^2(0,L)$ (i.e., $\Delta r = 0$ and $\Delta d = 0$). In the sequel, we denote by a subscript ``e'' the equilibrium value of the different quantities. For instance, $Z_e$ denotes the equilibrium value of $Z$. Noting that $u_{D,e} = u_e$ and $v_{D,e} = v_e$, we obtain
\begin{align*}
0 & = A_K Z_e + \Gamma_e \\
0 & = \lambda_j w_{j,e} + a_j u_e + b_j v_e + d_{e,j} , & j \geq n+1
\end{align*}
In particular, from $v_e = K Z_e$, we have 
\begin{equation*}
0 = A_K Z_e + \Gamma_e = A Z_e + e^{-D A} B v_e + \Gamma_e .
\end{equation*}
Since the first rows of $A$ and $\Gamma_e$ are null and $E_1 e^{-D A} B = 1$, we obtain $v_e = 0$ and
\begin{subequations}
\begin{align}
Z_e & = - A_K^{-1} \Gamma_e \label{eq: equilibrium - Ze} \\
u_e & = E_1 Z_e = - E_1 A_K^{-1} \Gamma_e \label{eq: equilibrium - ue} \\
w_{j,e} & = - \frac{a_j}{\lambda_j} u_e - \frac{d_{e,j}}{\lambda_j} , & j \geq n+1 \label{eq: equilibrium - wje}
\end{align}
We introduce 
$X_e = Z_e$
because $A X_e + B v_{D,e} + \Gamma_e = A_K Z_e + \Gamma_e = 0$, which is compatible with the Artstein transformation since $v_e = 0$ implies $Z_e = X_e + \int_{t-D}^t e^{(t-s-D)A} B v_e \,\mathrm{d}s$. The equilibrium condition of the integral component for reference tracking is given by 
\begin{equation}\label{eq: equilibrium - zetae}
\zeta_e = E_{n+2} X_e = - E_{n+2} A_K^{-1} \Gamma_e ,
\end{equation}
\end{subequations}
where $E_{n+2} = \begin{bmatrix} 0 & \ldots & 0 & 1 \end{bmatrix} \in \mathbb{R}^{1 \times (n+2)}$. Noting that $\lambda_j w_{j,e} = - a_j u_e - d_{e,j}$ for $j \geq n+1$ where $(a_j)_j$ and $(d_{e,j})_j$ are square-summable sequences and $\lambda_j \rightarrow + \infty$ when $j \rightarrow + \infty$, both $(w_{j,e})_j$ and $(\lambda_j w_{j,e})_j$ are square-summable sequences. Hence we define
\begin{equation}\label{eq: definition equilibrium we}
w_e \triangleq \sum\limits_{j \geq 1} w_{j,e} e_j \in D(\mathcal{A}) = H^2(0,L) \cap H^1_0(0,L) 
\end{equation}
which is convergent in $H_0^1(0,L)$. In particular, we obtain from the last line of $A X_e + \Gamma_e = 0$ and using (\ref{eq: equilibrium - wje}) that
\begin{align*}
& L_1 X_{1,e} = r_e + \sum\limits_{j \geq n+1} \frac{e_j'(0)}{\lambda_j} d_{e,j} \\
\Leftrightarrow\; & \frac{1}{L} u_e - \sum_{j \geq n+1} \frac{e_j'(0)}{\lambda_j} a_j u_e + \sum_{j=1}^{n} w_{j,e} e'_j(0) \\
& \hspace{4cm} = r_e + \sum\limits_{j \geq n+1} \frac{e_j'(0)}{\lambda_j} d_{e,j} \\
\Leftrightarrow\; &  \sum_{j \geq 1} w_{j,e} e'_j(0) + \frac{1}{L} u_e = r_e \\
\Leftrightarrow\; &  w'_e(0) + \frac{1}{L} u_e = r_e .
\end{align*}
Then, introducing $y_e \triangleq w_e + \frac{x}{L} u_e \in L^2(0,L)$, we obtain $y'_e(0) = r_e$, which corresponds to the desired reference tracking. Finally, since
\begin{multline*}
\mathcal{A} w_e = \sum\limits_{j \geq 1} \lambda_j w_{j,e} e_j 
= - \sum\limits_{j \geq 1} a_j e_j u_{e} - \sum\limits_{j \geq 1} d_{e,j} e_j  \\
= - a u_{e} - b v_{e} - d_e ,
\end{multline*}
we have $\mathcal{A} w_e + a u_{D,e} + b v_{D,e} + d_e = 0$.

\begin{remark}
The above developments show that the equilibrium point of the closed-loop infinite-dimensional system given by (\ref{eq: equilibrium - Ze}-\ref{eq: equilibrium - zetae}) and (\ref{eq: definition equilibrium we}) is fully determined by the constant values of the reference signal $r_e$ and the distributed disturbance $d_e$.
\end{remark}

\subsection{Dynamics of deviations}

We now define the deviations of the various quantities with respect to their equilibrium value: $\Delta X = X - X_e$, $\Delta Z = Z - Z_e$, $\Delta w = w - w_e$, $\Delta w_j = w_j - w_{j,e}$, $\Delta \zeta = \zeta - \zeta_e$, $\Delta u = u - u_e$ (first component of $\Delta Z$), $\Delta u_D = u_D - u_e$ (first component of $\Delta X$), $\Delta v = v - v_e$, and $\Delta v_D = v_D - v_{D,e}$. Then, in original coordinates:
\begin{equation}\label{eq: final model - PDE}
\Delta w_t = \mathcal{A} \Delta w + a \Delta u_D + b \Delta v_D + \Delta d
\end{equation}
and
\begin{align*}
\Delta\dot{X}(t) & = \mathcal{A} \Delta X(t) + B \Delta v_D(t) + \Delta\Gamma(t) \\
\Delta\dot{w}_j(t) & = \lambda_j \Delta w_j(t) + a_j \Delta u_D(t) + b_j \Delta v_D(t) + \Delta d_j(t)
\end{align*}
for $j \geq n+1$ with the auxiliary control input
$\Delta v(t) = \chi_{[0,+\infty)}(t) K \Delta Z(t)$
where 
\begin{equation}\label{eq: HL - Artstein transformation on Delta quantities}
\Delta Z(t) = \Delta X(t) + \int_{t-D}^t e^{(t-s-D)A} B \Delta v(s) \,\mathrm{d}s .
\end{equation}
In $Z$ coordinates, the closed-loop dynamics is given by
\begin{subequations}
\begin{align}
\Delta\dot{Z}(t) & = A_K \Delta Z(t) + \Delta\Gamma(t) \label{eq: final model - Z coordinates - Z} \\
\Delta\dot{w}_j(t) & = \lambda_j \Delta w_j(t) + a_j \Delta u_D(t) + b_j \Delta v_D(t) + \Delta d_j(t) \label{eq: final model - Z coordinates - w_j}
\end{align}
\end{subequations}
for $j \geq n+1$.

\section{Stability analysis}\label{sec: stability analysis}

\subsection{Main stability result}

The objective of this section is to establish the following stability result, taking the form of an Input-to-State Stability (ISS) estimate with fading memory  of both the reference input $r$ and the distributed perturbation $d$.

\begin{theorem}\label{thm: ISS closed-loop system}
There exist $\kappa, \overline{C}_1 > 0$ such that, for every $\epsilon \in [0,1)$, there exists $\overline{C}_2(\epsilon) > 0$ such that 
\begin{align}
& \Delta u_D(t)^2 + \Delta \zeta(t)^2 + \Vert \Delta w(t) \Vert_{H_0^1(0,L)}^2 \label{eq: HL - ISS estimate Delta_u_D Delta_zeta Delta_w} \\
& \qquad\leq 
\overline{C}_1 e^{- 2 \kappa t} \left( \Delta u_D(0)^2 + \Delta \zeta(0)^2 + \Vert \Delta w(0) \Vert_{H_0^1(0,L)}^2 \right) \nonumber \\
& \qquad\phantom{\leq}\,
+ \overline{C}_2(\epsilon) \sup\limits_{0 \leq s \leq t} e^{-2 \epsilon\kappa (t-s)} \{ \Delta r(s)^2 + \Vert \Delta d(s) \Vert^2 \} \nonumber .
\end{align}
Moreover, the constants $\kappa, \overline{C}_1, \overline{C}_2(\epsilon)$ can be chosen independently of $r_e$ and $d_e$.
\end{theorem}

Since $\Delta w(t,x) = \Delta y(t,x) - \frac{x}{L} \Delta u_D(t)$, we deduce from the continuous embedding $H_0^1(0,L) \subset L^\infty(0,L)$ (see, e.g., \cite{brezis2010functional}) 
the following corollary.
\begin{corollary}\label{cor: ISS closed-loop system - original coordinate}
Let $\kappa > 0$ be provided by Theorem~\ref{thm: ISS closed-loop system}. There exists $\tilde{C}_1 > 0$ such that, for every $\epsilon \in [0,1)$, there exists $\tilde{C}_2(\epsilon) > 0$ such that
\begin{align}
& \Vert \Delta y(t) \Vert_{L^\infty(0,L)} \label{eq: ISS closed-loop system - original coordinate} \\
& \quad\leq 
\tilde{C}_1 e^{- \kappa t} \left( \vert \Delta u_D(0) \vert + \vert \Delta \zeta(0) \vert + \Vert \Delta w(0) \Vert_{H_0^1(0,L)} \right) \nonumber \\
& \quad\phantom{\leq}\,
+ \tilde{C}_2(\epsilon) \sup\limits_{0 \leq s \leq t} e^{- \epsilon\kappa (t-s)} \{ \vert \Delta r(s) \vert + \Vert \Delta d(s) \Vert \} \nonumber .
\end{align}
\end{corollary}

We also deduce the following corollary concerning the asymptotic behavior of the closed-loop system in the case of convergent reference signal $r(t)$ and distributed disturbance $d(t)$ as $t \rightarrow + \infty$.

\begin{corollary}\label{cor: ISS closed-loop system}
Assume that $r(t) \rightarrow r_e$ and $d(t) \rightarrow d_e$ when $t \rightarrow + \infty$. Then $w(t) \rightarrow w_e$ in $H_0^1$ norm, $y(t) \rightarrow y_e$ in both $L^\infty$ and $L^2$ norm, $u(t) \rightarrow u_e$, and $\zeta(t) \rightarrow \zeta_e$ with exponential vanishing of the contribution of the initial conditions.
\end{corollary}

\begin{remark}
In the particular case $n = 0$, which corresponds to an exponentially stable open-loop reaction-diffusion equation (\ref{eq: investigated reaction-diffusion equation - beginning}-\ref{eq: investigated reaction-diffusion equation - ending}), the above results ensure that the stability of the closed-loop system is preserved after introduction of the two integral states $v$ and $z$. 
\end{remark}

In order to prove the claimed stability result, we resort as in~\cite{prieur2018feedback} to the Lyapunov function 
\begin{align}
V(t) 
& = \frac{M}{2} \Delta Z(t)^\top P \Delta Z(t) \label{eq: definition V} \\
& \phantom{=}\, + \frac{M}{2} \int_{\max(t-D,0)}^{t} \Delta Z(s)^\top P \Delta Z(s) \,\mathrm{d}s \nonumber \\
& \phantom{=}\, - \frac{1}{2} \sum_{j \geq 1} \lambda_j \Delta w_j(t)^2 , \nonumber
\end{align}
where $P \in \mathbb{R}^{(n+2) \times (n+2)}$ is the solution of the Lyapunov equation $A_K^\top P + P A_K = -I$ and $M > 0$ is chosen such that
$$
M > \max\left( \frac{\gamma_1 \lambda_1}{\lambda_m(P)} ,
4 \left( \gamma_1 \Vert a \Vert^2 + 2 \Vert b \Vert^2 \Vert e^{-D A_K} \Vert^2 \Vert K \Vert^2 \right) \right) 
$$
with $\gamma_1 \triangleq 2 \max\left( 1 , D e^{2 D \Vert A \Vert} \Vert B K \Vert^2 \right)$.

\begin{remark}
The first term in the definition (\ref{eq: definition V}) of $V$ accounts for the stability of the finite-dimensional truncated model (\ref{eq: final model - Z coordinates - Z}), expressed in $Z$ coordinates, capturing the $n$ first modes of the reaction-diffusion equation. The motivation behind the introduction of the second (integral) term relies on the fact that it allows, in conjunction with (\ref{eq: HL - Artstein transformation on Delta quantities}), the derivation of an upper-estimate of $\Vert \Delta X(t) \Vert$ (i.e., the state of the truncated model in its original $X$ coordinates) based on $V(t)$. Finally, the last term is used to capture the countable infinite number of modes of the original reaction-diffusion equation (\ref{eq: final model - PDE}), including those that where neglected in the control design. Note that $\langle \mathcal{A} \Delta w(t) , \Delta w(t) \rangle = \sum_{j \geq  1} \lambda_j \Delta w_j(t)^2$.
\end{remark}

\subsection{Preliminary Lemmas for the proof of Theorem \ref{thm: ISS closed-loop system}}

We derive hereafter various lemmas that will be useful in the sequel to establish the stability properties of the closed-loop system. First, we estimate $\Delta\Gamma(t)$ as follows.

\begin{lemma}
There exists a constant $M_d > 0$ such that
\begin{equation*}
\Vert \Delta\Gamma(t) \Vert^2 \leq M_d^2 ( \Delta r(t)^2 + \Vert \Delta d(t) \Vert^2 ) \quad\forall t\geq 0.
\end{equation*}
\end{lemma}

\textbf{Proof:}
By definition of $\Delta\Gamma(t)$ and using the Cauchy-Schwarz inequality we have
\begin{align*}
& \Vert \Delta\Gamma(t) \Vert^2 \\[-3mm]
& = \Vert \Delta D_1(t) \Vert^2 + \left\vert \Delta r(t) + \sum_{j \geq n+1} \frac{e_j'(0)}{\lambda_j} \Delta d_j(t) \right\vert^2 \\
& \leq \sum\limits_{j=1}^{n} \Delta d_j(t)^2
+ 2 \Delta r(t)^2
+ 2 \sum\limits_{j \geq n+1} \left\vert \frac{e_j'(0)}{\lambda_j} \right\vert^2 
\sum\limits_{j \geq n+1} \Delta d_j(t)^2 \\
& \leq M_d^2 ( \Delta r(t)^2 + \Vert \Delta d(t) \Vert^2 ) 
\end{align*}
with 
$
M_d^2 = {2\max\left( 1 , \sum\limits_{j \geq n+1} \left\vert \frac{e_j'(0)}{\lambda_j} \right\vert^2 \right)} < + \infty$,
since, by (\ref{eq: symptotic behaviors}), we have $\left\vert \frac{e_j'(0)}{\lambda_j} \right\vert^2 \sim \frac{2L}{\pi^2 j^2}$ when $j \rightarrow +\infty$.
\qed


\begin{lemma}\label{eq: upper estimate system trajectory}
There exists a constant $C_1 > 0$ such that
\begin{subequations}
\begin{align}
V(t) & \geq C_1 \sum_{j \geq 1} (1 + \vert \lambda_j \vert) \Delta w_j(t)^2 , \label{lem: upper estimate delta w_j with V} \\
V(t) & \geq C_1 \left( \Delta u_D(t)^2 + \Delta \zeta(t)^2 + \Vert \Delta w(t) \Vert_{H_0^1(0,L)}^2 \right) , \label{lem: upper estimate delta u delta zeta and delta w with V} \\
V(t) & \geq C_1 \Vert \Delta Z(t) \Vert^2 , \label{lem: upper estimate Delta Z with V}
\end{align}
\end{subequations}
for every $t \geq 0$.
\end{lemma}

\textbf{Proof:}
From (\ref{eq: HL - Artstein transformation on Delta quantities}) with $\Delta v = K \Delta Z$, we obtain that
\begin{align}
& \Vert \Delta X(t) \Vert^2 \nonumber \\
& \leq
2 \Vert \Delta Z(t) \Vert^2 \nonumber \\
& \phantom{\leq}\, + 2 D e^{2 D \Vert A \Vert} \Vert B K \Vert^2 \int_{\max(t-D,0)}^{t} \Vert \Delta Z(s) \Vert^2 \,\mathrm{d}s \nonumber \\
& \leq \gamma_1 \left( \Vert \Delta Z(t) \Vert^2 + \int_{\max(t-D,0)}^{t} \Vert \Delta Z(s) \Vert^2 \,\mathrm{d}s \right)  \label{eq: prel upper estimate norm Delta X1}
\end{align}
with $\gamma_1 = 2 \max\left( 1 , D e^{2 D \Vert A \Vert} \Vert B K \Vert^2 \right) > 0$. Thus, we have
\begin{align*}
& \Delta Z(t)^\top P \Delta Z(t) 
+ \int_{(t-D,t)\cap(0,+\infty)} \Delta Z(s)^\top P \Delta Z(s) \,\mathrm{d}s \\
& \hspace{5.5cm} \geq \frac{\lambda_m(P)}{\gamma_1} \Vert \Delta X(t) \Vert^2 .
\end{align*}
Noting that
\begin{align*}
\sum_{j \geq 1} \lambda_j \Delta w_j(t)^2
& \leq \sum_{j \geq n+1} \lambda_j \Delta w_j(t)^2 + \lambda_1 \sum_{j=1}^{n} \Delta w_j(t)^2 \\
& \leq \sum_{j \geq n+1} \lambda_j \Delta w_j(t)^2 + \lambda_1 \Vert \Delta X(t) \Vert^2 ,
\end{align*}
we obtain
\begin{equation*}
V(t) 
\geq \left( \frac{M \lambda_m(P)}{2\gamma_1} - \frac{\lambda_1}{2} \right) \Vert \Delta X(t) \Vert^2 
- \frac{1}{2} \sum_{j \geq n+1} \lambda_j \Delta w_j(t)^2 .
\end{equation*}
Since $M > \frac{\gamma_1 \lambda_1}{\lambda_m(P)} > 0$, we obtain the existence of $\gamma_2 = \frac{1}{2} \min\left( \frac{M \lambda_m(P)}{\gamma_1} - \lambda_1 , 1 \right) > 0$ such that
\begin{equation}\label{eq: intermediate lowe estimate V}
V(t) \geq \gamma_2 \Big( \Vert \Delta X(t) \Vert^2 - \sum_{j \geq n+1} \lambda_j \Delta w_j(t)^2 \Big),
\end{equation}
from which we obtain (\ref{lem: upper estimate delta w_j with V}). Now, as in~\cite{prieur2018feedback}, from the series expansions (\ref{eq: series expansion solutions}) and (\ref{eq: definition equilibrium we}) that are convergent in $H_0^1(0,L)$, we infer that
\begin{multline}\label{eq: computation norm H_0^1 of Delta w}
\Vert \Delta w(t) \Vert_{H_0^1(0,L)}^2 =
\sum\limits_{i,j \geq 1} \Delta w_i(t) \Delta w_j(t) \int_0^L e'_i(x) e'_j(x) \,\mathrm{d}x \\
= \int_0^L c(x) \Delta w(t,x)^2 \mathrm{d}x - \sum\limits_{j \geq 1} \lambda_j \Delta w_j(t)^2 , 
\end{multline}
where the second equality follows from an integration by part and the facts that $e''_j + c e_j = \lambda_j e_j$, $e_j(0)=e_j(L)=0$, and $(e_i)_{i \geq 1}$ is a Hilbert basis of $L^2(0,L)$. Hence, using the fact that $- \sum\limits_{1 \leq j \leq n} \lambda_j \Delta w_j(t)^2 \leq 0$, the following estimates hold:
\begin{align*}
& \Vert \Delta w(t) \Vert_{H_0^1(0,L)}^2 \\
& \leq \Vert c \Vert_{L^\infty(0,L)} \sum\limits_{j \geq 1} \Delta w_j(t)^2 - \sum\limits_{j \geq n+1} \lambda_j \Delta w_j(t)^2 \\
& \leq \Vert c \Vert_{L^\infty(0,L)} \sum\limits_{j = 1}^n \Delta w_j(t)^2
- \sum\limits_{j \geq n+1} \left( \lambda_j - \Vert c \Vert_{L^\infty(0,L)} \right) \Delta w_j(t)^2 \\
& \leq \gamma_3 \Big( \sum\limits_{j = 1}^n \Delta w_j(t)^2 - \sum\limits_{j \geq n+1} \lambda_j \Delta w_j(t)^2 \Big) 
\end{align*}
for some constant $\gamma_3 > 0$ because $\lambda_j \underset{j \rightarrow + \infty}{\longrightarrow} - \infty$ whence $ - \left( \lambda_j - \Vert c \Vert_{L^\infty(0,L)} \right) \sim - \lambda_j$ when $j \rightarrow + \infty$ with $\lambda_j < 0$ for all $j \geq n+1$. Therefore, we obtain from (\ref{eq: intermediate lowe estimate V}) that
 \begin{equation*}
V(t) 
\geq \gamma_2 \left( \Delta u_D(t)^2 + \Delta \zeta(t)^2 \right) 
+ \frac{\gamma_2}{\gamma_3}\Vert \Delta w(t) \Vert_{H_0^1(0,L)}^2 ,
\end{equation*}
which provides (\ref{lem: upper estimate delta u delta zeta and delta w with V}). Finally, from the definition of $V$ given by (\ref{eq: definition V}) and using (\ref{eq: prel upper estimate norm Delta X1}), we also have
\begin{align*}
V(t)
& \geq \frac{M \lambda_m(P)}{2} 
\left( \Vert \Delta Z(t) \Vert^2 + \int_{\max(t-D,0)}^{t} \Vert \Delta Z(s) \Vert^2 \,\mathrm{d}s \right) \\
& \phantom{=}\; \underbrace{- \frac{1}{2} \sum_{j \geq n+1} \lambda_j \Delta w_j(t)^2}_{\geq 0}
- \frac{\lambda_1}{2} \Vert \Delta X(t) \Vert^2 \\
& \geq \frac{1}{2} ( M \lambda_m(P) - \lambda_1 \gamma_1 ) \\
& \phantom{\geq}\; \times
\left( \Vert \Delta Z(t) \Vert^2 + \int_{\max(t-D,0)}^{t} \Vert \Delta Z(s) \Vert^2 \,\mathrm{d}s \right) \\
& \geq \gamma_1 \gamma_2 \Vert \Delta Z(t) \Vert^2 ,
\end{align*}
which gives (\ref{lem: upper estimate Delta Z with V}).
\qed

\subsection{End of proof of Theorem \ref{thm: ISS closed-loop system}}

We are now in a position to establish the stability properties of the closed-loop system and prove Theorem~\ref{thm: ISS closed-loop system}. We first study the exponential decay properties of $V$ for $t \geq D$.

\begin{lemma}\label{lem: exp decay V for t geq D}
There exist $\kappa, C_2 > 0$ such that, for every $\epsilon \in [0,1)$,
\begin{align*}
V(t) & \leq e^{-2 \kappa (t-D)} V(D) \\
& \phantom{\leq}\; + \frac{C_2}{1-\epsilon} \sup\limits_{0 \leq s \leq t} e^{-2\epsilon\kappa(t-s)} \left\{ \Delta r(s)^2 + \Vert \Delta d(s) \Vert^2 \right\} 
\end{align*}
for every $t \geq D$.
\end{lemma}

\textbf{Proof:}
First, we note that, for $t > D$, 
\begin{align*}
& \frac{\mathrm{d}}{\mathrm{d}t} \left[ \int_{t-D}^t \Delta Z(s)^\top P \Delta Z(s) \,\mathrm{d}s \right] \\
& \quad = \Delta Z(t)^\top P \Delta Z(t) - \Delta Z(t-D)^\top P \Delta Z(t-D) \\
& \quad = \int_{t-D}^t \frac{\mathrm{d}}{\mathrm{d}\tau} \left[ \Delta Z(\tau)^\top P \Delta Z(\tau) \right](s) \,\mathrm{d}s .
\end{align*}
Since $\mathcal{A}$ is self-adjoint and since $\langle \mathcal{A} \Delta w(t) , \Delta w(t) \rangle = \sum_{j \geq  1} \lambda_j \Delta w_j(t)^2$, we have for every $t > D$, 
\begin{align*}
& \dot{V}(t) \\
& = \frac{M}{2} \Delta Z(t)^\top \left( A_K^\top P + P A_K \right) \Delta Z(t) 
+ M \Delta Z(t)^\top P \Delta\Gamma(t) \\
& \phantom{=}\, + \frac{M}{2} \int_{t-D}^{t} \Delta Z(s)^\top \left( A_K^\top P + P A_K \right) \Delta Z(s) \,\mathrm{d}s  \\
& \phantom{=}\, + M \int_{t-D}^{t} \Delta Z(s)^\top P \Delta \Gamma(s) \,\mathrm{d}s 
- \langle \mathcal{A} \Delta w(t) , \Delta w_t(t) \rangle \\
& = - \frac{M}{2} \Vert \Delta Z(t) \Vert^2 
+ M \Delta Z(t)^\top P \Delta\Gamma(t) \\
& \phantom{=}\, - \frac{M}{2} \int_{t-D}^{t} \Vert \Delta Z(s) \Vert^2 \,\mathrm{d}s
+ M \int_{t-D}^{t} \Delta Z(s)^\top P \Delta \Gamma(s) \,\mathrm{d}s \\
& \phantom{=}\, - \Vert \mathcal{A} \Delta w(t) \Vert_{L^2(0,L)}^2 - \langle \mathcal{A} \Delta w(t) , a \rangle \Delta u_D(t)  \\
& \phantom{=}\, - \langle \mathcal{A} \Delta w(t) , b \rangle \Delta v_D(t) - \langle \mathcal{A} \Delta w(t) , \Delta d(t) \rangle .
\end{align*}
Using the following estimates:
\begin{multline*}
\Delta Z(t)^\top P \Delta\Gamma(t) 
\leq \Vert \Delta Z(t) \Vert \Vert P \Vert \Vert \Delta \Gamma(t) \Vert \\
\leq \frac{1}{4} \Vert \Delta Z(t) \Vert^2 + \Vert P \Vert^2 M_d^2 ( \Delta r(t)^2 + \Vert \Delta d(t) \Vert^2 ) ,
\end{multline*}
\begin{align*}
& \int_{t-D}^{t} \Delta Z(s)^\top P \Delta \Gamma(s) \,\mathrm{d}s \\
& \qquad\leq \int_{t-D}^{t} \Vert \Delta Z(s) \Vert \Vert P \Vert \Vert \Delta \Gamma(s)\Vert \,\mathrm{d}s \\
& \qquad\leq \frac{1}{4} \int_{t-D}^{t} \Vert \Delta Z(s) \Vert^2 \,\mathrm{d}s \\
& \qquad\phantom{\leq}\; + D \Vert P \Vert^2 M_d^2 \sup\limits_{t-D \leq s \leq t} \{ \Delta r(s)^2 + \Vert \Delta d(s) \Vert^2 \} ,
\end{align*}
\begin{align*}
& \left\vert \langle \mathcal{A} \Delta w(t) , a \rangle \Delta u_D(t) \right\vert \\
& \qquad\leq \frac{1}{4} \Vert \mathcal{A} \Delta w(t) \Vert_{L^2(0,L)}^2 + \Vert a \Vert^2 \vert \Delta u_D(t) \vert^2  \\
& \qquad\leq \frac{1}{4} \Vert \mathcal{A} \Delta w(t) \Vert_{L^2(0,L)}^2 + \Vert a \Vert^2 \Vert \Delta X(t) \Vert^2  \\
& \qquad\leq \frac{1}{4} \Vert \mathcal{A} \Delta w(t) \Vert_{L^2(0,L)}^2 \\
& \qquad\phantom{\leq}\, + \gamma_1 \Vert a \Vert^2 \left( \Vert \Delta Z(t) \Vert^2 + \int_{t-D}^t \Vert \Delta Z(s) \Vert^2 \,\mathrm{d}s \right)   ,
\end{align*}
\begin{align*}
& \left\vert \langle \mathcal{A} \Delta w(t) , b \rangle \Delta v_D (t) \right\vert \\
& \leq \frac{1}{4} \Vert \mathcal{A} \Delta w(t) \Vert_{L^2(0,L)}^2 + \Vert b \Vert^2 \vert \Delta v_D(t) \vert^2  \\
& \leq \frac{1}{4} \Vert \mathcal{A} \Delta w(t) \Vert_{L^2(0,L)}^2 + \Vert b \Vert^2 \Vert K \Vert^2 \Vert \Delta Z(t-D) \Vert^2  \\
& \leq \frac{1}{4} \Vert \mathcal{A} \Delta w(t) \Vert_{L^2(0,L)}^2 
+ 2 \Vert b \Vert^2 \Vert e^{-D A_K} \Vert^2 \Vert K \Vert^2 \Vert \Delta Z(t) \Vert^2  \\
& \phantom{\leq}\; + 2 M_d^2 D^2 e^{2D \Vert A_K \Vert} \Vert b \Vert^2 \Vert K \Vert^2 \sup\limits_{t-D \leq s \leq t} \{ \Delta r(s)^2 + \Vert \Delta d(s) \Vert^2 \} ,
\end{align*}
where we have used in the latter inequality that $\Delta\dot{Z}(t) = A_K \Delta Z(t) + \Delta \Gamma(t)$ whence 
\begin{equation*}
\Delta Z(t-D) = e^{- D A_K} \Delta Z(t) + \int_{t}^{t-D} e^{(t-D-s)A_K} \Delta\Gamma(s) \mathrm{d}s ,
\end{equation*}
and
\begin{align*}
\left\vert \langle \mathcal{A} \Delta w(t) , \Delta d(t) \rangle \right\vert 
& \leq \frac{1}{4} \Vert \mathcal{A} \Delta w(t) \Vert_{L^2(0,L)}^2 + \Vert \Delta d(t) \Vert^2 ,
\end{align*}
we obtain that, for every $t > D$,
\begin{multline*}
\dot{V}(t) \leq \left( - \frac{M}{4} + \gamma_1 \Vert a \Vert^2 + 2 \Vert b \Vert^2 \Vert e^{-D A_K} \Vert^2 \Vert K \Vert^2 \right) \\
\times \left( \Vert \Delta Z(t) \Vert^2 + \int_{t-D}^{t} \Vert \Delta Z(s) \Vert^2 \,\mathrm{d}s \right) \\
- \frac{1}{4} \Vert \mathcal{A} \Delta w(t) \Vert_{L^2(0,L)}^2 + \gamma_4 \sup\limits_{t-D \leq s \leq t} \{ \Delta r(s)^2 + \Vert \Delta d(s) \Vert^2 \} 
\end{multline*}
where 
\begin{equation*}
\gamma_4 = 1 + M_d^2 \left\{ (1+D) M \Vert P \Vert^2 + 2 D^2 e^{2 D \Vert A_K \Vert} \Vert b \Vert^2 \Vert K \Vert^2 \right\} .
\end{equation*}
Since $M > 4 \left( \gamma_1 \Vert a \Vert^2 + 2 \Vert b \Vert^2 \Vert e^{-D A_K} \Vert^2 \Vert K \Vert^2 \right)$, setting
\begin{equation*}
\gamma_5 = M/4 - \left( \gamma_1 \Vert a \Vert^2 + 2 \Vert b \Vert^2 \Vert e^{-D A_K} \Vert^2 \Vert K \Vert^2 \right) > 0
\end{equation*}
we have
\begin{align*}
& \dot{V}(t) \\
& \leq - \gamma_5 \left( \Vert \Delta Z(t) \Vert^2 + \int_{t-D}^{t} \Vert \Delta Z(s) \Vert^2 \,\mathrm{d}s \right) \\
& \phantom{\leq}\; - \frac{1}{4} \Vert \mathcal{A} \Delta w(t) \Vert_{L^2(0,L)}^2 + \gamma_4 \sup\limits_{t-D \leq s \leq t} \{ \Delta r(s)^2 + \Vert \Delta d(s) \Vert^2 \} \\ 
& \leq - \frac{\gamma_5}{\lambda_M(P)} \left( \Delta Z(t)^\top P \Delta Z(t)   
+ \int_{t-D}^{t} \Delta Z(s)^\top P \Delta Z(s) \,\mathrm{d}s \right) \\
& \phantom{\leq}\;  - \frac{1}{4} \Vert \mathcal{A} \Delta w(t) \Vert_{L^2(0,L)}^2 + \gamma_4 \sup\limits_{t-D \leq s \leq t} \{ \Delta r(s)^2 + \Vert \Delta d(s) \Vert^2 \} ,
\end{align*}
for every $t > D$. Now, since $\lambda_{j} \geq 0$ when $1 \leq j \leq n$ and $\lambda_j \leq \lambda_{n+1} < 0$ when $j \geq n+1$, we have, for every $t \geq 0$,
\begin{multline*}
- \sum_{j \geq 1} \lambda_j \Delta w_j(t)^2 
\leq - \sum_{j=n+1}^{+\infty} \lambda_j \Delta w_j(t)^2 \\
\leq \gamma_6 \sum_{j=n+1}^{+\infty} \lambda_j^2 \Delta w_j(t)^2
\leq \gamma_6 \Vert \mathcal{A} \Delta w(t) \Vert_{L^2(0,L)}^2 
\end{multline*}
with $\gamma_6 = 1/\vert \lambda_{n+1} \vert > 0$. Setting $\Delta p (s)^2 = \Delta r (s)^2 + \Vert \Delta d (s) \Vert^2$, we infer that
\begin{align*}
& \dot{V}(t) \\
& \leq - \frac{2 \gamma_5}{M \lambda_M(P)} \frac{M}{2} \\
& \phantom{\leq -}\; \times \left( \Delta Z(t)^\top P \Delta Z(t)   
+ \int_{t-D}^{t} \Delta Z(s)^\top P \Delta Z(s) \,\mathrm{d}s \right) \\
& \phantom{\leq}\; - \frac{1}{2} \frac{1}{2\gamma_6} \gamma_6 \Vert \mathcal{A} \Delta w(t) \Vert_{L^2(0,L)}^2 + \gamma_4 \sup\limits_{t-D \leq s \leq t} \Delta p (s)^2 \\
& \leq - 2 \kappa \frac{M}{2} \left( \Delta Z(t)^\top P \Delta Z(t)   
+ \int_{t-D}^{t} \Delta Z(s)^\top P \Delta Z(s) \,\mathrm{d}s \right) \\
& \phantom{\leq}\; - 2\kappa \frac{1}{2} \gamma_6 \Vert \mathcal{A} \Delta w(t) \Vert_{L^2(0,L)}^2 + \gamma_4 \sup\limits_{t-D \leq s \leq t} \Delta p (s)^2 \\ 
& \leq - 2 \kappa V(t) + \gamma_4 \sup\limits_{t-D \leq s \leq t} \{ \Delta r(s)^2 + \Vert \Delta d(s) \Vert^2 \} 
\end{align*}
for every $t > D$ where $\kappa = \frac{1}{2} \min\left( \frac{2 \gamma_5}{M \lambda_M(P)} , \frac{1}{2 \gamma_6} \right) > 0$. Then, we obtain, for every $t \geq D$ and every $\epsilon\in[0,1)$,
\begin{align*}
& V(t) - e^{-2\kappa(t-D)} V(D) \\
& \leq \gamma_4 e^{-2\kappa t} \int_{D}^{t} e^{2\kappa\tau} \sup\limits_{\tau-D \leq s \leq \tau} \Delta p (s)^2 \mathrm{d}\tau \\
& \leq \gamma_4 e^{-2\kappa t} \int_{D}^{t} e^{2(1-\epsilon)\kappa\tau} \mathrm{d}\tau\, \times \sup\limits_{D \leq \tau \leq t} \left[ e^{2\epsilon\kappa\tau} \sup\limits_{\tau-D \leq s \leq \tau} \Delta p (s)^2 \right] \\
& \leq \frac{\gamma_4}{2(1-\epsilon)\kappa} e^{-2\kappa t} e^{2(1-\epsilon)\kappa t}  \times \sup\limits_{D \leq \tau \leq t}\;\sup\limits_{\tau-D \leq s \leq \tau} e^{2\epsilon\kappa\tau} \Delta p (s)^2 \\
& \leq \frac{\gamma_4}{2(1-\epsilon)\kappa} e^{-2\epsilon\kappa t} \times \sup\limits_{D \leq \tau \leq t}\;\sup\limits_{\tau-D \leq s \leq \tau} e^{2\epsilon\kappa(s+D)} \Delta p (s)^2 \\
& \leq \frac{\gamma_4 e^{2\epsilon\kappa D}}{2(1-\epsilon)\kappa} e^{-2\epsilon\kappa t} \sup\limits_{D \leq \tau \leq t}\;\sup\limits_{\tau-D \leq s \leq \tau} e^{2\epsilon\kappa s} \Delta p (s)^2 \\
& \leq \frac{\gamma_4 e^{2\kappa D}}{2(1-\epsilon)\kappa} e^{-2\epsilon\kappa t} \sup\limits_{0 \leq s \leq t} e^{2\epsilon\kappa s} \{ \Delta r(s)^2 + \Vert \Delta d(s) \Vert^2 \} \\
& \leq \frac{\gamma_4 e^{2\kappa D}}{2(1-\epsilon)\kappa} \sup\limits_{0 \leq s \leq t} e^{-2\epsilon\kappa (t-s)} \{ \Delta r(s)^2 + \Vert \Delta d(s) \Vert^2 \} 
\end{align*}
where we have used to establish the fourth inequality that, for a given $\tau \in [D,t]$, $\tau - D \leq s \leq \tau$ implies $\tau \leq s+D$. The claimed estimate holds with $C_2 = \gamma_4 e^{2\kappa D} / (2\kappa)$.
\qed


\begin{lemma}\label{lem: upper estimate V for t leq D}
There exist constants $C_3,C_4 > 0$ such that
\begin{align*}
V(t) 
& \leq C_3 \left( \Delta u_D(0)^2 + \Delta \zeta(0)^2 + \Vert \Delta w(0) \Vert_{H_0^1(0,L)}^2 \right) \\
& \phantom{\leq}\; + C_4 \sup\limits_{0 \leq s \leq t} \{ \Delta r(s)^2 + \Vert \Delta d(s) \Vert^2 \} 
\end{align*}
for every $t\in[0,D]$ with $\Delta u_D(0) = - u_e$.
\end{lemma}

\textbf{Proof:}
For $0 \leq t \leq D$, we have 
\begin{align*}
V(t) 
& = \frac{M}{2} \left( \Delta Z(t)^\top P \Delta Z(t)   
+ \int_{0}^{t} \Delta Z(s)^\top P \Delta Z(s) \,\mathrm{d}s \right) \\
& \phantom{=}\; - \frac{1}{2}\sum_{j \geq 1} \lambda_j \Delta w_j(t)^2.
\end{align*}
We note that, for $0 \leq t < D$, $\Delta u_D(t) = u(t-D) - u_e = - u_e = \Delta u_D(0)$ and $\Delta v_D(t) = v(t-D) - v_e = 0$, whence
\begin{align*}
& \dot{V}(t) \\
& = \frac{M}{2} \Delta Z(t)^\top \left( A_K^\top P + P A_K \right) \Delta Z(t) 
+ M \Delta Z(t)^\top P \Delta \Gamma(t) \\
& \phantom{=}\; + \frac{M}{2} \Delta Z(t)^\top P \Delta Z(t) 
- \langle \mathcal{A} \Delta w(t) , \Delta w_t(t) \rangle \\
& = - \frac{M}{2} \Vert \Delta Z(t) \Vert^2 
+ M \Delta Z(t)^\top P \Delta \Gamma(t) \\
& \phantom{=}\; + \frac{M}{2} \Delta Z(t)^\top P \Delta Z(t)
- \Vert \mathcal{A} \Delta w(t) \Vert_{L^2(0,L)}^2 \\
& \phantom{=}\; - \langle \mathcal{A} \Delta w(t) , a \rangle \Delta u_D(t) - \langle \mathcal{A} \Delta w(t) , \Delta d(t) \rangle \\
& \leq \frac{M (\Vert P \Vert + \lambda_M(P) -1)}{2} \Vert \Delta Z(t) \Vert^2 
+ \frac{M \Vert P \Vert}{2} \Vert \Delta\Gamma(t) \Vert^2 \\
& \phantom{\leq}\; - \frac{1}{2} \Vert \mathcal{A} \Delta w(t) \Vert_{L^2(0,L)}^2 
+ \Vert a \Vert^2 \vert \Delta u_D(0) \vert^2
+ \Vert \Delta d(t) \Vert^2 \\
& \leq \frac{M (\Vert P \Vert + \lambda_M(P) -1)}{2} \Vert \Delta Z(t) \Vert^2 
+ \Vert a \Vert^2 \Vert \Delta X(0) \Vert^2 \\
& \phantom{\leq}\; + \max\left( 1 , \frac{M M_d^2 \Vert P \Vert}{2} \right) ( \Delta r(t)^2 + \Vert \Delta d(t) \Vert^2 ) .
\end{align*}
Noting that $\Delta\dot{Z}(t) = A_K \Delta Z(t) + \Delta\Gamma(t)$ and $\Delta Z(0) = \Delta X(0)$, we have 
\begin{equation*}
\Delta Z(t) = e^{ A_K t} \Delta X(0) + \int_0^t e^{ A_K (t-\tau)} \Delta \Gamma(\tau) \mathrm{d}\tau .
\end{equation*}
Thus, we obtain the existence of $\gamma_7, \gamma_8>0$ such that
\begin{equation*}
\dot{V}(t) \leq  \gamma_7 \Vert \Delta X(0) \Vert^2 + \gamma_8 \sup\limits_{0 \leq s \leq t} \{ \Delta r(s)^2 + \Vert \Delta d(s) \Vert^2 \} 
\end{equation*}
for $0 \leq t < D$. Therefore, 
\begin{multline*}
V(t)  \leq V(0) + D \gamma_7 \Vert \Delta X(0) \Vert^2 \\
+ D \gamma_8 \sup\limits_{0 \leq s \leq t} \{ \Delta r(s)^2 + \Vert \Delta d(s) \Vert^2 \} 
\end{multline*}
for  $0 \leq t \leq D$. To conclude, using (\ref{eq: computation norm H_0^1 of Delta w}), we estimate $V(0)+ D \gamma_7 \Vert \Delta X(0) \Vert^2$ as follows: 
\begin{align*}
& V(0) + D \gamma_7 \Vert \Delta X(0) \Vert^2 \\
& = \frac{M}{2} \Delta X(0)^\top P \Delta X(0) + D \gamma_7 \Vert \Delta X(0) \Vert^2 - \frac{1}{2} \sum_{j \geq 1} \lambda_j \Delta w_j(0)^2 \\
& \leq \frac{M \lambda_M(P)+2D \gamma_7}{2} \left( \Delta u_D(0)^2 + \Delta \zeta(0)^2 + \sum\limits_{j = 1}^{n} \Delta w_j(0)^2 \right) \\
& \phantom{\leq}\;  + \dfrac{1}{2} \left( \Vert \Delta w(0) \Vert_{H_0^1(0,L)}^2 + \Vert c \Vert_{L^\infty(0,L)} \Vert \Delta w(0) \Vert_{L^2(0,L)}^2 \right) \\
& \leq \frac{M \lambda_M(P)+2D \gamma_7}{2} \left( \Delta u_D(0)^2 + \Delta \zeta(0)^2 \right) \\
& \phantom{\leq}\; + \dfrac{1}{2} \left( 1 + L^2 \left\{ \Vert c \Vert_{L^\infty(0,L)} + M \lambda_M(P) + 2D \gamma_7 \right\} \right) \\
& \phantom{\leq\; +}\; \times \Vert \Delta w(0) \Vert_{H_0^1(0,L)}^2 ,
\end{align*}
where we have used the Poincar{\'e} inequality to derive the last estimate: $\Vert f \Vert_{L^2(0,L)} \leq L \Vert f \Vert_{H_0^1(0,L)}$ for every $f \in H_0^1(0,L)$. Combining the two latter estimates, the result follows.
\qed


\begin{lemma}\label{lem: exp decay V}
There exist $\kappa, C_5 > 0$ such that, for every $\epsilon \in [0,1)$, there exists $C_6(\epsilon) > 0$ such that
\begin{multline*}
V(t) \leq C_5 e^{-2 \kappa t} \left( \Delta u_D(0)^2 + \Delta \zeta(0)^2 + \Vert \Delta w(0) \Vert_{H_0^1(0,L)}^2 \right) \\
+ C_6(\epsilon) \sup\limits_{0 \leq s \leq t} e^{-2 \epsilon\kappa (t-s)} \{ \Delta r(s)^2 + \Vert \Delta d(s) \Vert^2 \} \quad \forall t\geq 0.
\end{multline*}
\end{lemma}

\textbf{Proof:}
When $0 \leq t \leq D$, Lemma~\ref{lem: upper estimate V for t leq D} yields
\begin{align*}
& V(t) \\
& \leq C_3 e^{2 \kappa D} e^{-2\kappa t} \left( \Delta u_D(0)^2 + \Delta \zeta(0)^2 + \Vert \Delta w(0) \Vert_{H_0^1(0,L)}^2 \right) \\
& \phantom{\leq}\;  + C_4 e^{2 \epsilon \kappa D} \sup\limits_{0 \leq s \leq t} e^{-2\epsilon\kappa(t-s)} \{ \Delta r(s)^2 + \Vert \Delta d(s) \Vert^2 \} 
\end{align*}
because $D-t \geq 0$ and $D-t+s \geq 0$ for all $0 \leq s \leq t \leq D$. When $t \geq D$, we infer from Lemma~\ref{lem: exp decay V for t geq D}, from the latter estimate evaluated in $t = D$, and by using again the notation $\Delta p (s)^2 = \Delta r (s)^2 + \Vert \Delta d (s) \Vert^2$, that
\begin{align*}
& V(t) \\
& \leq e^{-2 \kappa (t-D)} V(D) 
+ \frac{C_2}{1-\epsilon} \sup\limits_{0 \leq s \leq t} e^{-2\epsilon\kappa(t-s)} \Delta p (s)^2 \\
& \leq C_3 e^{-2 \kappa (t-D)} \left( \Delta u_D(0)^2 + \Delta \zeta(0)^2 + \Vert \Delta w(0) \Vert_{H_0^1(0,L)}^2 \right) \\
& \phantom{\leq}\; + C_4 e^{-2 \epsilon\kappa (t-D)} \sup\limits_{0 \leq s \leq D} e^{2\epsilon\kappa s} \{ \Delta r(s)^2 + \Vert \Delta d(s) \Vert^2 \} \\
& \phantom{\leq}\; + \frac{C_2}{1-\epsilon} \sup\limits_{0 \leq s \leq t} e^{-2\epsilon\kappa(t-s)} \left\{ \Delta r(s)^2 + \Vert \Delta d(s) \Vert^2 \right\} \\
& \leq C_3 e^{2 \kappa D} e^{-2 \kappa t} \left( \Delta u_D(0)^2 + \Delta \zeta(0)^2 + \Vert \Delta w(0) \Vert_{H_0^1(0,L)}^2 \right) \\
& \phantom{\leq}\; + \left( C_4 e^{2\epsilon\kappa D} + \frac{C_2}{1-\epsilon} \right) \sup\limits_{0 \leq s \leq t} e^{-2\epsilon\kappa(t-s)} \Delta p (s)^2 .
\end{align*}
The claimed estimate holds with $C_5 = C_3 e^{2 \kappa D}$ and $C_6(\epsilon) = C_4 e^{2\epsilon\kappa D} + \frac{C_2}{1-\epsilon}$.
\qed

We are now in a position to prove the main result of this section, namely the stability result stated in Theorem~\ref{thm: ISS closed-loop system}. Indeed, from Lemmas~\ref{eq: upper estimate system trajectory} and~\ref{lem: exp decay V}, we infer the existence of constants $\overline{C}_1 = C_5 / C_1 > 0$ and $\overline{C}_2(\epsilon) = C_6(\epsilon) / C_1 > 0$ such that (\ref{eq: HL - ISS estimate Delta_u_D Delta_zeta Delta_w}) holds. Similarly, we obtain the following estimates which will be useful in the next section concerning the tracking performance:
\begin{align}
& \sum_{j \geq 1} (1 + \vert \lambda_j \vert) \Delta w_j(t)^2 \label{eq: HL - ISS estimate sum lambda_i w_i^2} \\
& \leq \overline{C}_1 e^{- 2 \kappa t} \left( \Delta u_D(0)^2 + \Delta \zeta(0)^2  + \Vert \Delta w(0) \Vert_{H_0^1(0,L)}^2 \right) \nonumber \\
& \phantom{\leq}\; + \overline{C}_2(\epsilon) \sup\limits_{0 \leq s \leq t} e^{-2 \epsilon\kappa (t-s)} \{ \Delta r(s)^2 + \Vert \Delta d(s) \Vert^2 \} , \nonumber
\end{align}
for every $t \geq 0$, and, as $\Delta v(t) = K \Delta Z(t)$ for $t \geq 0$ and $\Delta v(t) = 0$ for $t<0$,
\begin{align}
& \Vert \Delta v_D(t) \Vert^2 \label{eq: HL - ISS estimate Delta alpha} \\
& \leq \hat{C}_1 e^{- 2 \kappa t} \left( \Delta u_D(0)^2 + \Delta \zeta(0)^2 + \Vert \Delta w(0) \Vert_{H_0^1(0,L)}^2 \right) \nonumber \\
& \phantom{\leq}\; + \hat{C}_2(\epsilon) \sup\limits_{0 \leq s \leq t} e^{-2 \epsilon\kappa (t-s)} \{ \Delta r(s)^2 + \Vert \Delta d(s) \Vert^2 \} , \nonumber
\end{align}
for every $t \geq 0$ with $\hat{C}_1 = \Vert K \Vert^2 \overline{C}_1 e^{2 \kappa D}$ and $\hat{C}_2(\epsilon) = \Vert K \Vert^2 \overline{C}_2(\epsilon) e^{2 \epsilon\kappa D}$. This concludes the proof of Theorem \ref{thm: ISS closed-loop system}.

\begin{remark}
All the constants $C_i$ for $1 \leq i \leq 6$, and thus $\overline{C}_1,\overline{C}_2(\epsilon)$, defined in this section are independent of the considered equilibrium condition characterized by the quantities $r_e$ and $d_e$. Consequently, one can apply the result of Theorem~\ref{thm: ISS closed-loop system}, for the same values of the constants $\overline{C}_1,\overline{C}_2(\epsilon)$, to distinct equilibrium points associated with different constant values $r_e$ and $d_e$ successively taken by the reference signal $r(t)$ and the distributed disturbance $d(t)$, respectively. This feature will be illustrated in numerical computations in Section~\ref{sec: numerical illustration}.
\end{remark}

\section{Setpoint reference tracking analysis}\label{sec: setpoint reference tracking analysis}

It now remains to assess that the setpoint tracking of the reference signal $r(t)$ is achieved in the presence of the distributed disturbance $d(t)$. Specifically, we establish in this section the following tracking result.
\begin{theorem}\label{thm: ISS regulation of the output}
Let $\kappa > 0$ be provided by Theorem~\ref{thm: ISS closed-loop system}. There exists $\overline{C}_3 > 0$ such that, for every $\epsilon \in [0,1)$, there exists $\overline{C}_4(\epsilon) > 0$ such that 
\begin{align}
& \vert y_x(t,0) - r(t) \vert \label{eq: ISS estimate tracking performance} \\
& \leq 
\overline{C}_3 e^{- \kappa t} \left( \vert \Delta u_D(0) \vert + \vert \Delta \zeta(0) \vert + \Vert \Delta w(0) \Vert_{H_0^1(0,L)} \right. \nonumber \\
& \phantom{\leq}\, \hspace{5cm} \left. + \Vert \mathcal{A} \Delta w(0) \Vert_{L^2(0,L)} \right) \nonumber \\
& \phantom{\leq}\,
+ \overline{C}_4(\epsilon) \sup\limits_{0 \leq s \leq t} e^{- \epsilon\kappa (t-s)} \{ \vert \Delta r(s) \vert + \Vert \Delta d(s) \Vert + \Vert \dot{d}(s) \Vert \} \nonumber .
\end{align}
Moreover, the constants $\overline{C}_3, \overline{C}_4(\epsilon)$ can be chosen independently of the parameters $r_e$ and $d_e$.
\end{theorem}


\begin{corollary}\label{cor: ISS regulation of the output}
Assume that $r(t) \rightarrow r_e$, $d(t) \rightarrow d_e$, and $\dot{d}(t) \rightarrow 0$ when $t \rightarrow + \infty$. Then $y_x(t,0) \rightarrow r_e$ with exponential vanishing of the contribution of the initial conditions.
\end{corollary}

\begin{remark}
In the particular case $n = 0$, which corresponds to an exponentially stable open-loop reaction-diffusion equation (\ref{eq: investigated reaction-diffusion equation - beginning}-\ref{eq: investigated reaction-diffusion equation - ending}), the above results also ensure that the proposed control strategy achieves the setpoint reference tracking of the reference signal $r(t)$ while preserving the stability of the closed-loop system.
\end{remark}


\textbf{Proof of Theorem~\ref{thm: ISS regulation of the output}.}
Based on the identity $w_{e,x}(0) + \frac{1}{L} u_e = r_e$, we have the estimates:
\begin{multline}\label{eq: estimate y_x(t,0) - y_e_x(t,0)}
\left\vert y_x(t,0) - r(t) \right\vert
\quad\leq
\left\vert w_x(t,0) + \frac{1}{L}u_D(t) - r_e \right\vert
+ \vert \Delta r(t) \vert  \\
\leq 
\vert w_x(t,0) - w_{e,x}(0) \vert
+ \frac{1}{L} \vert \Delta u_D (t) \vert 
+ \vert \Delta r(t) \vert . 
\end{multline}
From the estimate of $\Delta u_D(t)$ provided by (\ref{eq: HL - ISS estimate Delta_u_D Delta_zeta Delta_w}), it is sufficient to study the term $w_x(t,0) - w_{e,x}(0) = \sum\limits_{j \geq 1} \Delta w_j(t) e_j'(0)$. Since $e_j'(0) \sim \sqrt{2/L} \sqrt{\vert \lambda_j \vert}$, there exists a constant $\gamma_9 > 0$ such that $\vert e'_j(0) \vert \leq \gamma_9 \sqrt{\vert \lambda_j \vert}$ for all $j \geq n + 1$. Let $m \geq n+1$ be such that $\eta \triangleq - \lambda_m > \kappa > 0$. Thus $\lambda_j \leq - \eta < - \kappa < 0$ for all $j \geq m$. We infer from the Cauchy-Schwarz inequality that
\begin{multline}\label{eq: estimate w_x(t,0) - w_e_x(t,0)}
\vert w_x(t,0) - w_{e,x}(0) \vert 
\leq \sum\limits_{j \geq 1} \vert \Delta w_j(t) \vert \vert e_j'(0) \vert  \\[-2mm]
\leq \sum\limits_{j = 1}^{m-1} \vert \Delta w_j(t) \vert \vert e_j'(0) \vert + \gamma_9 \sum\limits_{j \geq m} \sqrt{\vert \lambda_j \vert} \vert \Delta w_j(t) \vert \\
\leq \sqrt{\sum\limits_{j = 1}^{m-1} e_j'(0)^2} \sqrt{\sum\limits_{j = 1}^{m-1} \Delta w_j(t)^2}  \\
+ \gamma_9 \sqrt{\sum\limits_{j \geq m} \frac{1}{\vert \lambda_j \vert}} \sqrt{\sum\limits_{j \geq m} \lambda_j^2 \Delta w_j(t)^2} 
\end{multline}
where $\sum\limits_{j \geq m} \frac{1}{\vert \lambda_j \vert} < + \infty$ because $\lambda_j \sim -\pi^2 j^2 / L^2$. Based on (\ref{eq: HL - ISS estimate sum lambda_i w_i^2}), it is sufficient to study the term $\sqrt{\sum\limits_{j \geq m} \lambda_j^2 \Delta w_j(t)^2}$. To do so, we integrate for $j \geq m$ the dynamics (\ref{eq: final model - Z coordinates - w_j}) of the coefficient $\Delta w_j(t)$ as follows: 
\begin{align}
& \lambda_j \Delta w_j(t) = e^{\lambda_j t} \lambda_j \Delta w_j(0) \label{eq: integrated form of lambda_j Delta w_j} \\
&\quad + \int_{0}^{t} \lambda_j e^{\lambda_j(t-\tau)} \{ a_j \Delta u_D(\tau) + b_j \Delta v_D(\tau) + \Delta d_j(\tau) \} \,\mathrm{d}\tau . \nonumber
\end{align}
Now, integrating by parts and noting that $\Delta \dot{d}_j(\tau) = \dot{d}_j(\tau)$, we have
\begin{align}
& \int_{0}^{t} \lambda_j e^{\lambda_j(t-\tau)} \Delta d_j(\tau) \,\mathrm{d}\tau \label{eq: integral d_j term - IPP} \\
& \quad = - \Delta d_j(t) + e^{\lambda_j t} \Delta d_j(0) + \int_{0}^{t} e^{\lambda_j(t-\tau)} \dot{d}_j(\tau) \,\mathrm{d}\tau , \nonumber
\end{align} 
whence, 
\begin{align*}
& \vert \lambda_j \Delta w_j(t) \vert \\
& \leq e^{\lambda_j t} \vert \lambda_j \Delta w_j(0) \vert \\ 
& \phantom{\leq}\; + \int_0^t (-\lambda_j) e^{\lambda_j(t-\tau)} \left\{ \vert a_j \vert \vert \Delta u_D(\tau) \vert + \vert b_j \vert \vert \Delta v_D(\tau) \vert \right\} \,\mathrm{d}\tau \\
& \phantom{\leq}\; + \vert \Delta d_j(t) \vert + e^{\lambda_j t} \vert \Delta d_j(0) \vert + \int_{0}^{t} e^{\lambda_j(t-\tau)} \vert \dot{d}_j(\tau) \vert \,\mathrm{d}\tau \\
& \leq e^{-\eta t} \vert \lambda_j \Delta w_j(0) \vert 
+ \vert a_j \vert \int_0^t (-\lambda_j) e^{\lambda_j(t-\tau)} \vert \Delta u_D(\tau) \vert \,\mathrm{d}\tau \\
& \phantom{\leq}\; + \vert b_j \vert \int_0^t (-\lambda_j) e^{\lambda_j(t-\tau)} \vert \Delta v_D(\tau) \vert \,\mathrm{d}\tau \\
& \phantom{\leq}\; + \vert \Delta d_j(t) \vert + e^{-\eta t} \vert \Delta d_j(0) \vert + \int_{0}^{t} e^{-\eta(t-\tau)} \vert \dot{d}_j(\tau) \vert \,\mathrm{d}\tau .
\end{align*}
Now, as $\lambda_j \leq - \eta < - \kappa < - \epsilon\kappa$, the use of estimate (\ref{eq: HL - ISS estimate Delta_u_D Delta_zeta Delta_w}) and the introduction of the notations $\Delta\mathrm{CI} = \sqrt{\Delta u_D(0)^2 + \Delta \zeta(0)^2 + \Vert \Delta w(0) \Vert_{H_0^1(0,L)}^2}$ and $\Delta p(s) = \sqrt{ \Delta r(s)^2 + \Vert \Delta d(s) \Vert^2 }$ yield
\begin{align*}
& \int_0^t (-\lambda_j) e^{\lambda_j(t-\tau)} \vert \Delta u_D(\tau) \vert \,\mathrm{d}\tau \\
& \leq  (-\lambda_j) \sqrt{\overline{C}_1} e^{\lambda_j t} \int_{0}^{t} e^{- \lambda_j \tau} e^{-\kappa\tau} \,\mathrm{d}\tau\, \Delta\mathrm{CI} \\
& \phantom{\leq}\; + (-\lambda_j) \sqrt{\overline{C}_2(\epsilon)} e^{\lambda_j t} \int_{0}^{t} e^{- \lambda_j \tau} \sup\limits_{0 \leq s \leq \tau} e^{-\epsilon\kappa (\tau-s)} \Delta p(s) \,\mathrm{d}\tau \\
& \leq  (-\lambda_j) \sqrt{\overline{C}_1} e^{\lambda_j t} \int_{0}^{t} e^{- (\lambda_j+\kappa)\tau} \,\mathrm{d}\tau \,\Delta\mathrm{CI} \\
& \phantom{\leq}\; + (-\lambda_j) \sqrt{\overline{C}_2(\epsilon)} e^{\lambda_j t} \int_{0}^{t} e^{- (\lambda_j + \epsilon\kappa) \tau} \sup\limits_{0 \leq s \leq \tau} e^{\epsilon\kappa s} \Delta p(s) \,\mathrm{d}\tau \\
& \leq  \frac{\lambda_j}{\lambda_j + \kappa} \sqrt{\overline{C}_1} e^{\lambda_j t} (e^{-(\lambda_j + \kappa)t}-1) \,\Delta\mathrm{CI} \\
& \phantom{\leq}\; + \frac{\lambda_j}{\lambda_j + \epsilon\kappa} \sqrt{\overline{C}_2(\epsilon)} e^{\lambda_j t} (e^{-(\lambda_j + \epsilon\kappa)t}-1) \sup\limits_{0 \leq s \leq t} e^{\epsilon\kappa s} \Delta p(s) \\
& \leq  \frac{\eta}{\eta - \kappa} \sqrt{\overline{C}_1} e^{- \kappa t} \Delta\mathrm{CI} \\
& \phantom{\leq}\; + \frac{\eta}{\eta - \epsilon\kappa} \sqrt{\overline{C}_2(\epsilon)} e^{- \epsilon\kappa t} \sup\limits_{0 \leq s \leq t} e^{\epsilon\kappa s} \Delta p(s)\\
& \leq  \frac{\eta}{\eta - \kappa} \sqrt{\overline{C}_1} e^{- \kappa t} \Delta\mathrm{CI} \\
& \phantom{\leq}\; + \frac{\eta}{\eta - \epsilon\kappa} \sqrt{\overline{C}_2(\epsilon)} \sup\limits_{0 \leq s \leq t} e^{-\epsilon\kappa (t-s)} \Delta p(s) .
\end{align*}
Similarly, the use of the estimate (\ref{eq: HL - ISS estimate Delta alpha}) yields
\begin{align*}
& \int_0^t (-\lambda_j) e^{\lambda_j(t-\tau)} \vert \Delta v_D(\tau) \vert \,\mathrm{d}\tau \\
& \leq  \frac{\eta}{\eta - \kappa} \sqrt{\hat{C}_1} e^{- \kappa t} \Delta\mathrm{CI} \\
& \phantom{\leq}\; + \frac{\eta}{\eta - \epsilon\kappa} \sqrt{\hat{C}_2(\epsilon)} \sup\limits_{0 \leq s \leq t} e^{-\epsilon\kappa (t-s)} \Delta p(s) .
\end{align*}
Finally, since $\eta > \kappa$, we also infer from the Cauchy-Schwarz inequality that
\begin{align*}
& \int_{0}^{t} e^{-\eta(t-\tau)} \vert \dot{d}_j(\tau) \vert \,\mathrm{d}\tau = \int_{0}^{t} e^{-(\eta-\kappa)(t-\tau)} e^{-\kappa(t-\tau)} \vert \dot{d}_j(\tau) \vert \,\mathrm{d}\tau \\
& \leq \sqrt{\int_{0}^{t} e^{-2(\eta-\kappa)(t-\tau)} \,\mathrm{d}\tau} \sqrt{\int_{0}^{t} e^{-2\kappa(t-\tau)} \vert \dot{d}_j(\tau) \vert^2 \,\mathrm{d}\tau} \\
& \leq \sqrt{\frac{1}{2(\eta-\kappa)}} \sqrt{\int_{0}^{t} e^{-2\kappa(t-\tau)} \vert \dot{d}_j(\tau) \vert^2 \,\mathrm{d}\tau}.
\end{align*}
We infer from the three above estimates that
\begin{align*}
& \vert \lambda_j \Delta w_j(t) \vert \\
& \leq e^{-\eta t} \vert \lambda_j \Delta w_j(0) \vert \\
& \phantom{\leq} + \frac{\eta}{\eta-\kappa} \left(\vert a_j \vert \sqrt{\overline{C}_1} + \vert b_j \vert \sqrt{\hat{C}_1}\right) e^{-\kappa t} \Delta\mathrm{CI} \\
& \phantom{\leq} + \frac{\eta}{\eta-\epsilon\kappa} \left( \vert a_j \vert \sqrt{\overline{C}_2(\epsilon)} + \vert b_j \vert \sqrt{\hat{C}_2(\epsilon)} \right) \sup\limits_{0 \leq s \leq t} e^{-\epsilon\kappa (t-s)} \Delta p (s) \\
& \phantom{\leq} + \vert \Delta d_j(t) \vert + e^{-\eta t} \vert \Delta d_j(0) \vert \\
& \phantom{\leq} + \sqrt{\frac{1}{2(\eta-\kappa)}} \sqrt{\int_{0}^{t} e^{-2\kappa(t-\tau)} \vert \dot{d}_j(\tau) \vert^2 \mathrm{d}\tau}.
\end{align*}
Consequently we have:
\begin{align*}
& \vert \lambda_j \Delta w_j(t) \vert^2 \\
& \leq 6 e^{-2\eta t} \vert \lambda_j \Delta w_j(0) \vert^2 \\ 
& \phantom{\leq}\; + \frac{6 \eta^2}{(\eta-\kappa)^2} \left(\vert a_j \vert \sqrt{\overline{C}_1} + \vert b_j \vert \sqrt{\hat{C}_1}\right)^2 e^{-2\kappa t} \Delta\mathrm{CI}^2 \\
& \phantom{\leq}\; + \frac{6 \eta^2}{(\eta-\epsilon\kappa)^2} \left( \vert a_j \vert \sqrt{\overline{C}_2(\epsilon)} + \vert b_j \vert \sqrt{\hat{C}_2(\epsilon)} \right)^2 \\ 
& \phantom{\leq}\; \phantom{+}\; \times \sup\limits_{0 \leq s \leq t} e^{-2\epsilon\kappa (t-s)} \Delta p(s)^2 \\
& \phantom{\leq}\; + 6 \vert \Delta d_j(t) \vert^2 + 6 e^{-2\eta t} \vert \Delta d_j(0) \vert^2 \\
& \phantom{\leq}\; + \frac{3}{\eta-\kappa} \int_{0}^{t} e^{-2\kappa(t-\tau)} \vert \dot{d}_j(\tau) \vert^2 \mathrm{d}\tau 
\end{align*}
whence
\begin{align*}
& \sum\limits_{j \geq m} \lambda_j^2 \Delta w_j(t)^2 \\
& \leq 6 e^{- 2 \kappa t} \Vert \mathcal{A} \Delta w(0) \Vert_{L^2(0,L)}^2 \nonumber \\
& \ + \frac{12\eta^2}{(\eta-\epsilon\kappa)^2} (\Vert a \Vert^2 \overline{C}_1 + \Vert b \Vert^2 \hat{C}_1) e^{-2\kappa t} \Delta\mathrm{CI}^2 \nonumber \\
& \ + \frac{12 \eta^2}{(\eta-\epsilon\kappa)^2} (\Vert a \Vert^2 \overline{C}_2(\epsilon) + \Vert b \Vert^2 \hat{C}_2(\epsilon)) \sup\limits_{0 \leq s \leq t} e^{-2\epsilon\kappa (t-s)} \Delta p(s)^2 \nonumber \\
& \ + 6 \Vert \Delta d(t) \Vert^2 + 6 e^{-2 \kappa t} \Vert \Delta d(0) \Vert^2 \\
& \ + \frac{3}{\eta-\kappa} \int_{0}^{t} e^{-2\kappa(t-\tau)} \Vert \dot{d}(\tau) \Vert^2 \,\mathrm{d}\tau \nonumber \\
& \leq 6 e^{- 2 \kappa t} \Vert \mathcal{A} \Delta w(0) \Vert_{L^2(0,L)}^2 \nonumber \\ 
& \ + \frac{12 \eta^2}{(\eta-\epsilon\kappa)^2} (\Vert a \Vert^2 \overline{C}_1 + \Vert b \Vert^2 \hat{C}_1) e^{-2\kappa t} \Delta\mathrm{CI}^2 \nonumber \\
& \ + \frac{12 \eta^2}{(\eta-\epsilon\kappa)^2} (\Vert a \Vert^2 \overline{C}_2(\epsilon) + \Vert b \Vert^2 \hat{C}_2(\epsilon)) \sup\limits_{0 \leq s \leq t} e^{-2\epsilon\kappa (t-s)} \Delta p(s)^2 \nonumber \\
& \ + 12 \sup\limits_{0 \leq s \leq t} e^{-2\kappa(t-\tau)} \Vert \Delta d(\tau) \Vert^2 \\
& \ + \frac{3}{2(1-\epsilon)(\eta-\kappa)\kappa} \sup\limits_{0 \leq s \leq t} e^{-2\epsilon\kappa(t-s)} \Vert \dot{d}(s) \Vert^2  \nonumber 
\end{align*}
where we have used that
\begin{align*}
& \int_{0}^{t} e^{-2\kappa(t-\tau)} \Vert \dot{d}(\tau) \Vert^2 \,\mathrm{d}\tau \\
& \quad = \int_{0}^{t} e^{-2(1-\epsilon)\kappa(t-\tau)} e^{-2\epsilon\kappa(t-\tau)} \Vert \dot{d}(\tau) \Vert^2 \,\mathrm{d}\tau \\
& \quad \leq \int_{0}^{t} e^{-2(1-\epsilon)\kappa(t-\tau)} \,\mathrm{d}\tau \times \sup\limits_{0 \leq s \leq t} e^{-2\epsilon\kappa(t-s)} \Vert \dot{d}(s) \Vert^2 \\
& \quad \leq \frac{1}{2(1-\epsilon)\kappa} \sup\limits_{0 \leq s \leq t} e^{-2\epsilon\kappa(t-s)} \Vert \dot{d}(s) \Vert^2 .
\end{align*}
We deduce the existence of constants $C_7,C_8(\epsilon) > 0$ such that
\begin{align}
& \sum\limits_{j \geq m} \lambda_j^2 \Delta w_j(t)^2 \label{eq: estimate sum lambda_j^2 Delta w_j^2} \\
& \leq C_7 e^{-2\kappa t} \left( \Delta u_D(0)^2 + \Delta \zeta(0)^2 + \Vert \Delta w(0) \Vert_{H_0^1(0,L)}^2 \right. \nonumber \\
& \phantom{\leq}\; \hspace{5cm} \left. + \Vert \mathcal{A} \Delta w(0) \Vert_{L^2(0,L)}^2 \right) \nonumber \\
& \  + C_8(\epsilon) \sup\limits_{0 \leq s \leq t} e^{-2\epsilon\kappa (t-s)} \left( \Delta r(s)^2 + \Vert \Delta d(s) \Vert^2 + \Vert \dot{d}(s) \Vert^2 \right) . \nonumber  
\end{align}
Using now (\ref{eq: estimate y_x(t,0) - y_e_x(t,0)}) along with (\ref{eq: estimate w_x(t,0) - w_e_x(t,0)}) and estimates (\ref{eq: HL - ISS estimate Delta_u_D Delta_zeta Delta_w}), (\ref{eq: HL - ISS estimate sum lambda_i w_i^2}), and (\ref{eq: estimate sum lambda_j^2 Delta w_j^2}), we obtain the existence of the claimed constants $C_3,C_4(\epsilon) > 0$ such that the estimate (\ref{eq: ISS estimate tracking performance}) holds.  
\qed

\begin{remark}
At first sight, it might seem surprising that the estimate (\ref{eq: ISS estimate tracking performance}) on the tracking performance only involves the time derivative $\dot{d}$ of the distributed disturbance but not the time derivative $\dot{r}$ of the reference signal. Such a dissimilarity between the reference signal and the distributed disturbance is due to the explicit occurrence of the distributed perturbation $d$ in the dynamics (\ref{eq: final model - Z coordinates - w_j}) of the coefficient of projection $\Delta w_j$. Indeed, in order to estimate the term $\vert \lambda_j \Delta w_j(t) \vert$ from (\ref{eq: integrated form of lambda_j Delta w_j}), one needs to estimate the term $\int_{0}^{t} \lambda_j e^{\lambda_j(t-\tau)} \Delta d_j(\tau) \mathrm{d}\tau$. To do so, one first needs to eliminate the multiplicative factor $\lambda_j$ using, e.g., either an integration or an integration by parts. Simultaneously, we need to use Parseval identity in order to gather all coefficients $\Delta d_j(t)$. However, contrarily to the constant coefficients $a_j,b_j$, each coefficient $\Delta  d_j(t)$ is a function of time and thus cannot be pulled out of the integral. This remark motivates the integration by parts carried out in (\ref{eq: integral d_j term - IPP}). This way, the multiplicative factor $\lambda_j$ is eliminated and the subsequent estimates can be obtained. However, this is at the price of the emergence of the term $\dot{d}$ in the resulting tracking estimate. 
\end{remark}

\section{Numerical illustration}\label{sec: numerical illustration}

We take $c = 1.25$, $L = 2 \pi$, and $D = 1\,\mathrm{s}$. The three first eigenvalues of the open-loop system are $\lambda_1 = 1$, $\lambda_2 = 0.25$, and $\lambda_2 = -1$. Only the two first modes need to be stabilized. Thus we have $n = 2$ and we compute the feedback gain $K \in \mathbb{R}^{1 \times 4}$ such that the poles of the closed-loop truncated model (capturing the two unstable modes of the infinite-dimensional system plus two integral components, one for the control input and one for the reference tracking) are given by $-0.5$, $-0.6$, $-0.7$, and $-0.8$. The adopted numerical scheme is the modal approximation of the infinite-dimensional system using its first 10 modes. The initial condition is $y_0(x) = \frac{x}{L} \left( 1 - \frac{x}{L} \right)$.

\begin{figure}
     \centering
     	\subfigure[State $y(t,x)$]{
		\includegraphics[width=3.5in]{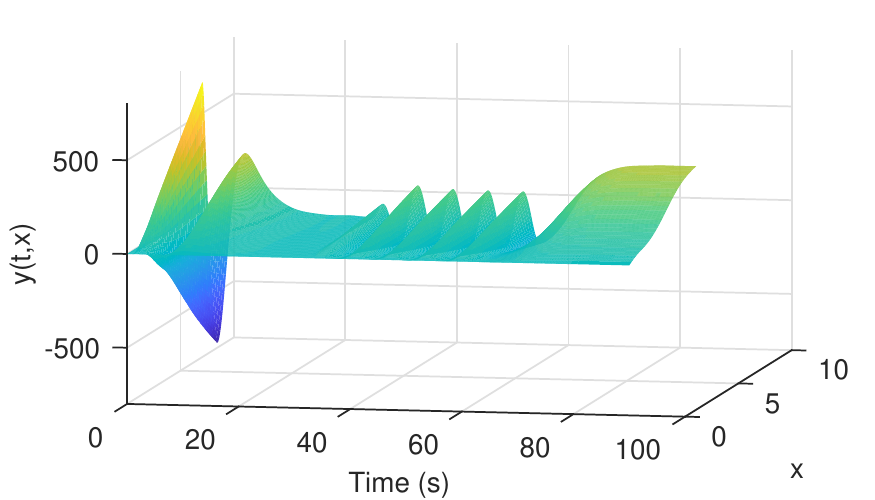}
		\label{fig: sim state X - time varying ref}
		}
     	\subfigure[Output $y_x(t,0)$ tracking the time-varying reference signal $r(t)$]{
		\includegraphics[width=3.5in]{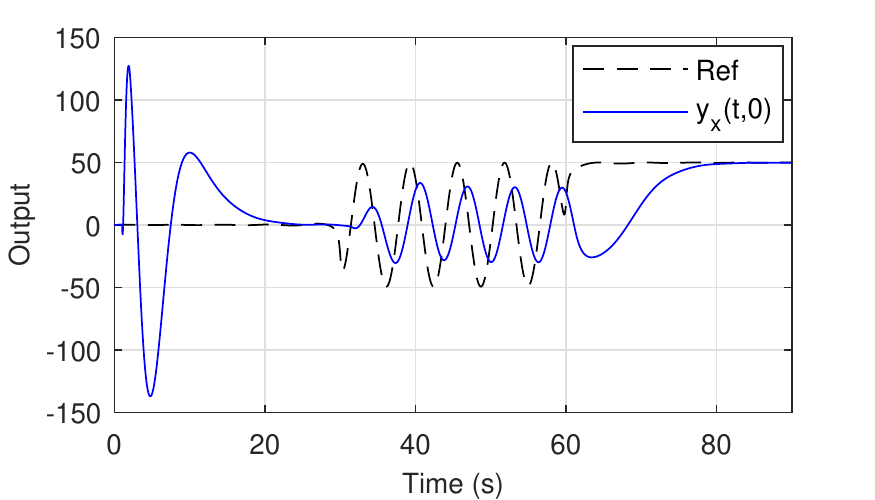}
		\label{fig: sim output y - time varying ref}
		}
     	\subfigure[Delayed control input $u(t-D)$]{
		\includegraphics[width=3.5in]{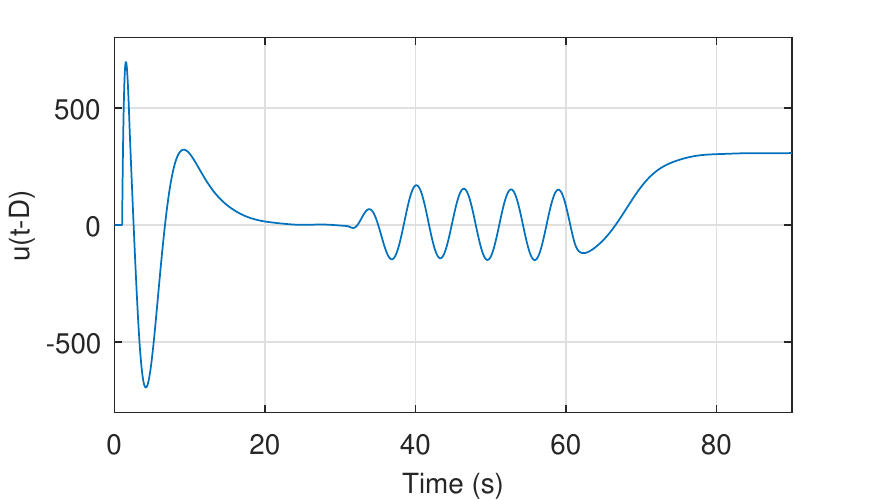}
		\label{fig: sim input u - time varying ref}
		}
     \caption{Time evolution of the closed-loop system for a time-varying reference signal $r(t)$ and a constant distributed perturbation $d(t,x) = x$}
     \label{fig: sim closed-loop - time varying ref}
\end{figure}

The simulation results for a time-varying reference $r(t)$ evolving within the range $[0,50]$ and the constant distributed disturbance $d(t,x) = x$ are depicted in Fig.~\ref{fig: sim closed-loop - time varying ref}. Applying first the obtained stability results for $t < 30\,\mathrm{s}$ with $r_e = 0$ and $d_e(x) = x$, we obtain that $y \rightarrow 0$ in $L^\infty(0,L)$ norm, $u(t) \rightarrow 0$, and $y_x(t,0) \rightarrow 0$, when $t \rightarrow + \infty$. This is compliant with the simulation result observed for increasing values of $t$ approaching $t = 30\,\mathrm{s}$. Consequently, the numerical simulation confirms that the proposed control strategy achieves the exponential stabilization of the closed-loop system while ensuring a zero steady-state left Neumann trace. Then, for $30\,\mathrm{s} < t < 60\,\mathrm{s}$, the tracking error remains bounded in the presence of an oscillatory reference signal. Finally, for $t > 60\,\mathrm{s}$, we apply again the obtained stability results but for $r_e = 50$ and $d_e(x) = x$. This time, we obtain that $y \rightarrow y_e \neq 0$ in $L^\infty(0,L)$ norm and $u(t) \rightarrow u_e \neq 0$ as $y_x(t,0) \rightarrow r_e = 50$ when $t \rightarrow + \infty$. In particular, conforming to the obtained ISS estimates with fading memory (\ref{eq: HL - ISS estimate Delta_u_D Delta_zeta Delta_w}-\ref{eq: ISS closed-loop system - original coordinate}) and (\ref{eq: ISS estimate tracking performance}), the impact of the variations of the reference signal around its nominal value $r_e$, i.e., configuration for which $\Delta r(t) \neq 0$, are eliminated as $t$ increases due to the action of the PI controller. This result provides a numerical confirmation of the efficiency of the proposed PI control strategy for the regulation control of the left Neumann trace of the system.

\begin{figure}
     \centering
     	\subfigure[State $y(t,x)$]{
		\includegraphics[width=3.5in]{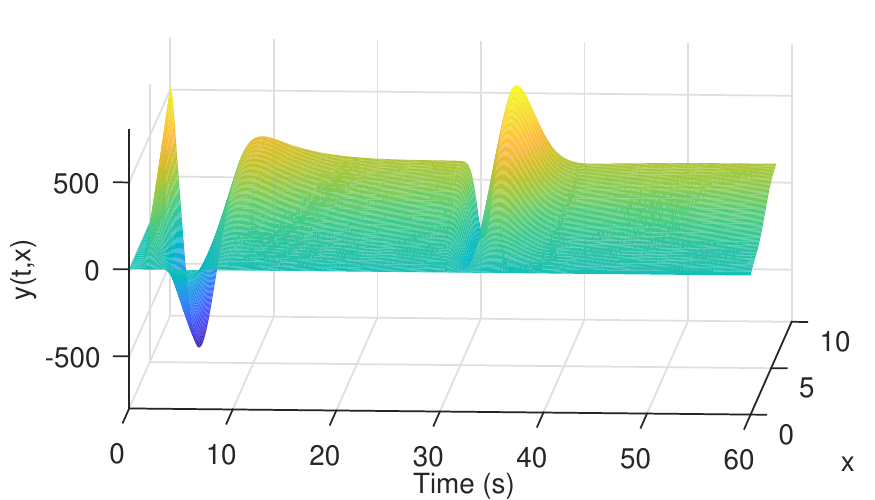}
		\label{fig: sim state X - time varying pert}
		}
     	\subfigure[Output $y_x(t,0)$ tracking the constant reference signal $r(t) = 50$]{
		\includegraphics[width=3.5in]{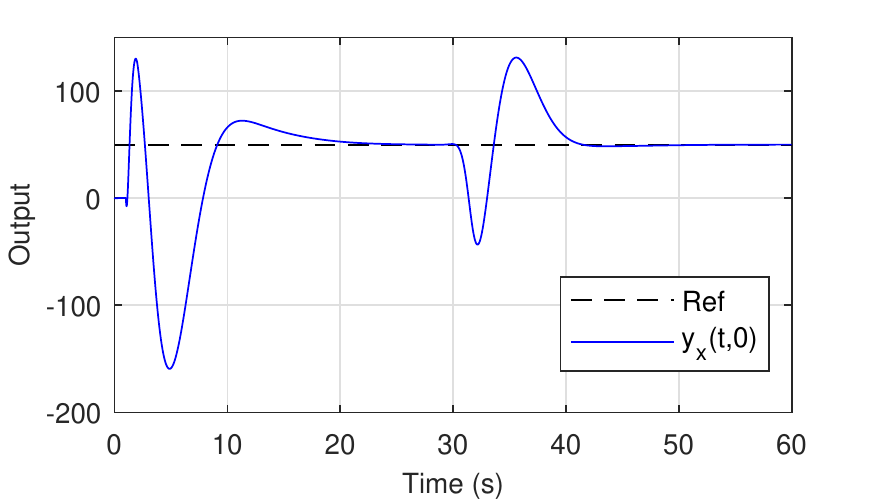}
		\label{fig: sim output y - time varying pert}
		}
     	\subfigure[Delayed control input $u(t-D)$]{
		\includegraphics[width=3.5in]{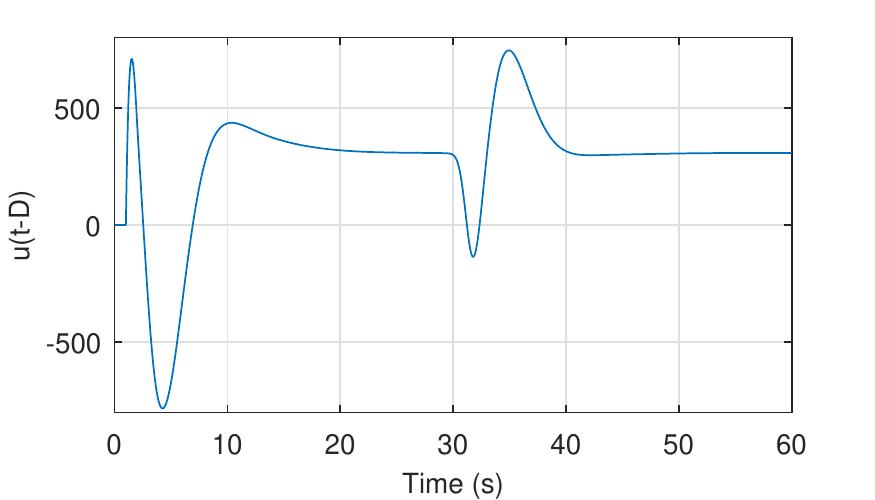}
		\label{fig: sim input u - time varying pert}
		}
     	\subfigure[Time varying component $d_0(t)$ of the distributed disturbance $d(t) = d_0(t) x$]{
		\includegraphics[width=3.5in]{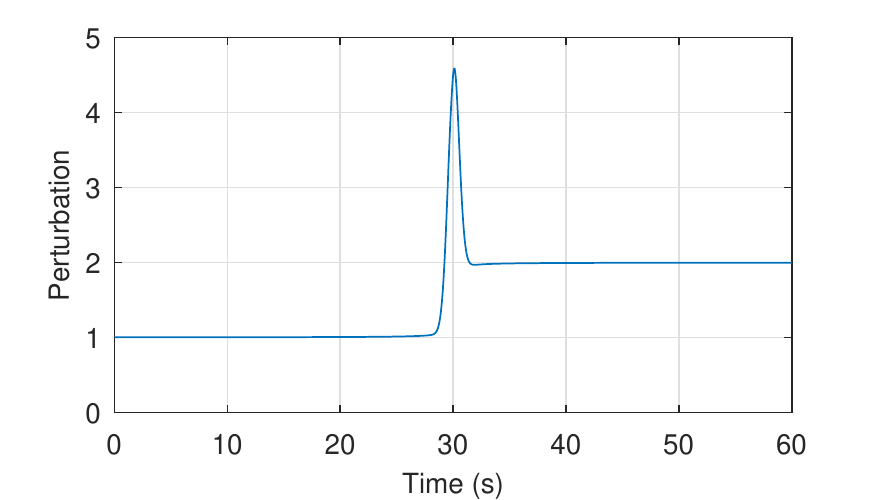}
		\label{fig: sim dist d - time varying pert}
		}
     \caption{Time evolution of the closed-loop system for a constant reference signal $r(t) = 50$ and a time-varying distributed disturbance $d(t,x) = d_0(t) x$}
     \label{fig: sim closed-loop - time varying pert}
\end{figure}

The simulation results for a constant reference $r(t)=50$ and the time-varying distributed disturbance $d(t,x) = d_0(t) x$ with $d_0$ given by Fig.~\ref{fig: sim dist d - time varying pert} are depicted in Fig.~\ref{fig: sim closed-loop - time varying pert}.  First, for time $t < 30\,\mathrm{s}$, the left Neumann trace successfully tracks the reference signal $r(t) = r_e = 50$ under the constant distributed perturbation $d(t,x) = d_e(x) = x$. Then, around $t = 30\,\mathrm{s}$, the magnitude of the perturbation increases from $d_0(t) = 1$ to $d_0(t) = 2$ with an overshoot around the value of $4.5$. This time-varying perturbation induces a perturbation on the setpoint reference tracking of the the left Neumann trace over the time-range $[30,40]\,\mathrm{s}$. However, once the perturbation reaches its steady-state value $d(t,x) = d_e(x) = 2x$, the impact of off-equilibrium perturbations are eliminated providing $y_x(t,0) \rightarrow r_e = 50$. This is compliant with the obtained ISS estimates with fading memory (\ref{eq: HL - ISS estimate Delta_u_D Delta_zeta Delta_w}-\ref{eq: ISS closed-loop system - original coordinate}) and (\ref{eq: ISS estimate tracking performance}) as the contribution of the variations of the perturbation around its nominal value $d_e$, i.e., configuration for which $\Delta d(t) \neq 0$, are eliminated as $t$ increases due to the action of the PI controller.

\section{Conclusion}\label{sec: conclusion}
We have achieved the PI regulation control of the left Neumann trace of a one dimensional linear reaction-diffusion equation with delayed right Dirichlet boundary control. The proposed control design approach extends to PI control a recently proposed approach for the delay boundary feedback control of infinite-dimensional systems via spectral reduction. Specifically, a finite-dimensional model capturing the unstable modes of the open-loop system has been obtained by spectral decomposition. Based on the classical Artstein transformation (used to handle the delay in the control input) and the pole schifting theorem, a PI controller has been derived. Then, the stability of the full infinite-dimensional closed-loop system has been assessed by using an adequate Lyapunov function, yielding an exponential Input-to-State Stability (ISS) estimate with fading memory of the time-varying reference signal and the time-varying distributed disturbance. Finally, a similar exponential ISS estimate with fading memory has been derived for the setpoint regulation control of the reference signal by the left Neumann trace.

As a conclusion, we indicate here potential directions for the extension of the work reported in this paper. 

First, it would be of interest to investigate whether the proposed PI control strategy can be used for the delay boundary regulation control of analogous PDEs. Good candidates in this direction are the linear Kuramoto-Sivashinsky equation \cite{guzman2019stabilization} and the wave equation as studied in~\cite{coron2006global}. 

Second, the work presented here was devoted to the control of a 1-D reaction diffusion. A natural research direction relies in the investigation of whether the proposed PI boundary control strategy could be applied to a multi-dimensional reaction-diffusion equation. This is not a straightforward extension of the developments presented in this paper since, in particular, we instrumentally used the fact that, in 1-D, $\sum 1/\vert \lambda_j \vert < + \infty$. Such a condition fails in multi-D. 

Finally, we assumed in this work that the measure of the full state is available. Future developments may be concerned with the development of an observer and the study of the stability of the resulting closed-loop system.


%
%
%

\appendix[Technical lemma]
The following lemma generalizes the result of~\cite[Chap.~12.4]{khalil2002nonlinear} to the case $D \neq 0$.

\begin{lemma}\label{lem: technical lemma}
Let $A \in \R^{n \times n}$, $B \in \R^{n \times m}$, $C  \in \R^{p \times n}$, and $D \in \R^{p \times m}$ be given matrices. The two following properties are equivalent:
\begin{itemize}
\item[(i)] The pair $(A,B)$ satisfies the Kalman condition and $\mathrm{rank}\begin{pmatrix} A & B \\ C & D \end{pmatrix} = n+p$.
\item[(ii)] The pair $\left( \begin{pmatrix} A & 0_{n\times p} \\ C & 0_{p\times p} \end{pmatrix} , \begin{pmatrix} B \\ D \end{pmatrix} \right)$ satisfies the Kalman condition.
\end{itemize}
\end{lemma}

In order to prove Lemma~\ref{lem: technical lemma}, we will use the following result.

\begin{lemma}\label{lem: technical lemma - prel}
Let $M \in \R^{q \times q}$ and $N \in \R^{q \times r}$ be given matrices. Assume that $(M,N)$ satisfies the Kalman condition. Then $\mathrm{Ran}\begin{pmatrix} M & N \end{pmatrix}=\R^q$, i.e., the matrix $\begin{pmatrix} M & N \end{pmatrix}$ is surjective.
\end{lemma}

\textbf{Proof of Lemma~\ref{lem: technical lemma - prel}:}
Noting that the surjectivity of $\begin{pmatrix} M & N \end{pmatrix}$ is equivalent to the condition $\ker M^\top \cap \ker N^\top = \{0\}$, let $\psi\in\R^n$ be such that $\psi^\top \begin{pmatrix} M & N \end{pmatrix}=0$. We have then $\psi^\top N=0$ and $\psi^\top M=0$, hence $\psi^\top M^k=0$ for every $k\in\N^*$ and thus $\psi M^k N=0$ for every $k\in\N$. Since $(M,N)$ satisfies the Kalman condition, we infer that $\psi=0$.
\qed


\textbf{Proof of Lemma~\ref{lem: technical lemma}:}
($(i) \Rightarrow (ii)$) Let us prove that, if $\psi_1\in\R^n$ et $\psi_2\in\R^p$ are such that
\begin{equation*}
\begin{pmatrix} \psi_1^\top & \psi_2^\top\end{pmatrix} \begin{pmatrix} B & AB & A^2B & \cdots & A^{n+p-1}B \\ D & CB & CAB & \cdots & CA^{n+p-2}B \end{pmatrix} = 0 ,
\end{equation*}
then $\psi_1=\psi_2=0$. Indeed, we have then $\psi_1^\top B+\psi_2^\top D=0$ and $(\psi_1^\top A+\psi_2^\top C)A^{i}B = 0$ for all $0 \leq i \leq n+p-2$. Noting that $n+p-2 \geq n-1$ and since $(A,B)$ satisfies the Kalman condition, we obtain that $\psi_1^\top A+\psi_2^\top C=0$. Consequently
\begin{equation*}
\begin{pmatrix} \psi_1^\top & \psi_2^\top\end{pmatrix}\begin{pmatrix} A & B \\ C & D \end{pmatrix} = 0 
\end{equation*} 
and hence $\psi_1=\psi_2=0$ since $\begin{pmatrix} A & B \\ C & D \end{pmatrix}$ is surjective.

\smallskip

\noindent($(ii) \Rightarrow (i)$) Let us prove that if $\psi\in\R^n$ is such that $\psi^\top \begin{pmatrix} B & AB & \cdots & A^{n-1}B\end{pmatrix}=0$ then $\psi=0$. Indeed, we first infer from the Hamilton-Cayley theorem that $\psi^\top A^kB=0$ for every $k\in\N$, and then
\begin{equation*}
\begin{pmatrix} \psi^\top & 0\end{pmatrix} \begin{pmatrix} B & AB & A^2B & \cdots & A^{n+p-1}B \\ D & CB & CAB & \cdots & CA^{n+p-2}B \end{pmatrix} = 0 .
\end{equation*}
Since the pair $\left( \begin{pmatrix} A & 0 \\ C & 0 \end{pmatrix} , \begin{pmatrix} B \\ D \end{pmatrix} \right)$ satisfies the Kalman condition, we obtain $\psi=0$. This shows that $(A,B)$ satisfies the Kalman condition. Besides, by Lemma~\ref{lem: technical lemma - prel}, since $\left( \begin{pmatrix} A & 0 \\ C & 0 \end{pmatrix} , \begin{pmatrix} B \\ D \end{pmatrix} \right)$ satisfies the Kalman condition, $\begin{pmatrix} A & 0 & B \\ C & 0 & D \end{pmatrix}$ is surjective. The lemma is proved.
\qed


\ifCLASSOPTIONcaptionsoff
  \newpage
\fi



\bibliographystyle{IEEEtranS}
\nocite{*}
\bibliography{IEEEabrv,mybibfile}

\end{document}